\documentclass[a4paper,centertags,oneside,11pt]{amsart}
\usepackage{mathrsfs}
\usepackage{appendix}
\usepackage{amssymb}
\usepackage{fancyhdr}
\usepackage{charter}
\usepackage{typearea}
\usepackage{pdfsync}
\usepackage{mathrsfs}
\usepackage[a4paper,top=3cm,bottom=3cm,left=3cm,right=3cm]{geometry}

\usepackage{enumitem}
\setlist[1]{itemsep=5pt}
\newcommand{\comment}[1]{}

\makeatletter
      \def\@setcopyright{}
      \def\serieslogo@{}
      \makeatother

\newcommand{\Complex}{\mathbb C}
\newcommand{\Real}{\mathbb R}

\newcommand{\ddbar}{\overline\partial}
\newcommand{\pr}{\partial}
\newcommand{\ol}{\overline}
\newcommand{\Td}{\widetilde}
\newcommand{\norm}[1]{\left\Vert#1\right\Vert}
\newcommand{\abs}[1]{\left\vert#1\right\vert}

\newcommand{\set}[1]{\left\{#1\right\}}
\newcommand{\To}{\rightarrow}

\newtheorem{theorem}{Theorem}[section]
\newtheorem{lemma}[theorem]{Lemma}
\newtheorem{corollary}[theorem]{Corollary}
\newtheorem{proposition}[theorem]{Proposition}
\newtheorem{question}[theorem]{Question}
\newtheorem{definition}[theorem]{Definition}

\newtheorem{example}[theorem]{Example}
\newtheorem{remark}[theorem]{Remark}
\numberwithin{equation}{section}
\begin{document}
\title[]{Szeg\H{o} kernel asymptotics and Morse inequalities on CR manifolds with $S^1$ action}
\author[]{Chin-Yu Hsiao}
\address{Institute of Mathematics, Academia Sinica and National Center for Theoretical Sciences, Astronomy-Mathematics Building, No. 1, Sec. 4, Roosevelt Road, Taipei 10617, Taiwan}
\thanks{Chin-Yu Hsiao was partially supported by Taiwan Ministry of Science of Technology project 104-2628-M-001-003-MY2, the Golden-Jade fellowship of Kenda Foundation and Academia Sinica Career Development Award.}
\email{chsiao@math.sinica.edu.tw or chinyu.hsiao@gmail.com}
\author[]{Xiaoshan Li}
\address{School of Mathematics
and Statistics, Wuhan University, Hubei 430072, China \& Institute of Mathematics, Academia Sinica, 6F, Astronomy-Mathematics Building,
No.1, Sec.4, Roosevelt Road, Taipei 10617, Taiwan}
\thanks{Xiaoshan Li was supported by  National Natural Science Foundation of China (Grant No. 11501422)}
\email{xiaoshanli@whu.edu.cn or xiaoshanli@math.sinica.edu.tw}

\dedicatory{Dedicated to Professor Ngaiming Mok for his 60th birthday}

\begin{abstract}
Let $X$ be a compact  connected CR manifold of dimension $2n-1, n\geq 2$. We assume that there is a transversal CR locally free $S^1$ action on $X$. Let $L^k$ be the $k$-th  power of a rigid CR line bundle $L$ over $X$. Without any assumption on the Levi-form of $X$, we obtain a scaling upper-bound for the partial Szeg\H{o} kernel on $(0,q)$-forms with values in $L^k$. After integration, this gives the weak Morse inequalities. By a refined spectral analysis, we also obtain the strong Morse inequalities in CR setting. We apply the strong Morse inequalities to show that the Grauert-Riemenschneider criterion is also true in the CR setting.
\end{abstract}

\maketitle \tableofcontents

\section{Introduction and statement of the main results}

The problem of embedding CR manifolds is prominent in areas such as complex analysis, partial differential equations and differential geometry. Let $X$ be a compact CR manifold of dimension $2n-1$, $n\geq 2$. When $X$ is strongly pseudoconvex and dimension of $X$ is greater than or equal to five, a classical theorem of
L. Boutet de Monvel~\cite{BdM1:74b} asserts that $X$ can be globally CR embedded into $\Complex^N$, for some $N\in\mathbb N$.
For a strongly pseudoconvex CR  manifold of dimension greater than five, the dimension of the kernel of the tangential Cauchy-Riemmann operator $\ddbar_b$ is infinite and we can find many CR functions to embed $X$ into complex space.
Inspired by Kodaira, the first-named author and Marinescu introduced in \cite{HM12} the idea of embedding CR manifolds by means of CR sections of tensor powers $L^k$ of a CR line bundle $L\To X$. To study Kodaira type embedding theorems on CR manifolds, it is crucial to be able to know

\begin{question}\label{qI}
When ${\rm dim\,}H^0_b(X,L^k)\gtrsim k^n$, for $k$ large, where $H^0_b(X,L^k)$ denotes the space of global smooth CR sections of $L^k$.
\end{question}

Inspired by Demailly \cite{D85,D91} (see also Getzler \cite{Ge89}), the first-named author and Marinescu established in \cite{HM12} analogues of the holomorphic Morse inequalities of Demailly for CR manifolds.

\begin{theorem}[Theorem 1.8, \cite{HM12}] \label{tI}
We assume that $Y(0)$ and $Y(1)$ hold at each point of $X$.
Then as $k\to\infty$,
\begin{equation}\label{*}
\begin{split}
&-{\rm dim\,}H^0_b(X, L^k)+{\rm dim\,}H^1_b(X, L^k)\\
&\quad\leqslant  \frac{k^{n}}{(2\pi)^{n}}\Bigr(-\int_X\int_{\Real_{x,0}}\abs{\det(\mathcal R_x^L+2s\mathcal{L}_x)}ds\,dv_X(x)\\
&\quad+\int_X\int_{\Real_{x,1}}\abs{\det(\mathcal R_x^L+2s\mathcal{L}_x)}ds\,dv_X(x)\Bigr)+o(k^n),
\end{split}
\end{equation}
where $\mathcal R_x^L$ is the associated curvature of $L$ at $x\in X$, $H^1_b(X, L^k)$ denotes the first $\ddbar_b$ cohomology group with values in $L^k$, $\mathcal{L}_x$ denotes the Levi form of $X$ at $x\in X$, and for $x\in X$, $q=0,1$,
\begin{equation}\label{***}\begin{split}
&\Real_{x,q}=\{s\in\Real;\, \mbox{$\mathcal R_x^L+2s\mathcal{L}_x$ has exactly $q$ negative eigenvalues} \\
&\quad\mbox{and $n-1-q$ positive eigenvalues}\}.
\end{split}\end{equation}
\end{theorem}
When $Y(0)$ and $Y(1)$ hold, from Kohn's results we know that ${\rm dim\,}H^0_b(X,L^k)<\infty$ and ${\rm dim\,}H^1_b(X,L^k)<\infty$. From \eqref{*}, we see that if
\begin{equation}\label{**}
\int_X\int_{\Real_{x,0}}\abs{\det(\mathcal R_x^L+2s\mathcal{L}_x)}dsdv_X(x)>
\int_X\int_{\Real_{x,1}}\abs{\det(\mathcal R_x^L+2s\mathcal{L}_x)}dsdv_X(x)
\end{equation}
then $L$ is big, that is ${\rm dim\,}H^0_b(X,L^k)\gtrsim k^n$. This is a very general criterion and it is desirable to refine it in some cases where \eqref{**} is not easy to verify. In general, it is very difficult to see when \eqref{**} holds even $L$ is positive.
The problem comes from the presence of positive eigenvalues of $\mathcal R^L_x$ and negative eigenvalues of $\mathcal{L}_x$.
By using Theorem~\ref{tI} to approach Question~\ref{qI}, we always have to impose extra conditions linking the Levi form and the curvature of the line bundle $L$. Similar problems also appear in the works of Marinescu~\cite{M96, M16}, Berman~\cite{Be05} where they studied the $\ddbar$-Neumann cohomology groups associated to a high power of a given holomorphic line bundle on a compact complex manifold with boundary. In order to get many holomorphic sections, they also have to
assume that, close to the boundary, the curvature of the line bundle is adapted to the Levi form of the boundary. In~\cite{Hsiao15}, by carefully studying semi-classical behaviour of microlocal Fourier transforms of the extreme functions for the spaces of lower energy forms of the associated Kohn Laplacian, the first-named author prove that $L$ is big when $L$ is positive,  $Y(0)$ and $Y(1)$ hold on $X$ under certain Sasakian conditions on $X$ and $L$ without any extra condition linking the Levi form of $X$ and the curvature of $L$. All these developments need the assumptions that the Levi form satisfies condition $Y(0)$ and $Y(1)$.

However, in some important problems in CR geometry, we need to know when $L$ is big without any assumption of the Levi form. For example, Ohsawa and Sibony~\cite{OS00} studied Kodaira type embedding theorems on Levi-flat CR manifolds. In their work, it is important to understand the space $H^0_b(X,L^k)$ for $k$ large. Adachi~\cite{Ad13} constructed a positive CR line bundle $L$ over a Levi-flat compact CR manifold $X$ of dimension $2n-1$ such that ${\rm dim\,}H^0_b(X,L^k)\lesssim k^{n-1}<k^n$ for $k$ large. We are lead to ask

\begin{question}\label{qII}
Can we establish some kind of Morse inequalities and Grauert-Riemenschneider criterion on some class of CR manifolds without any Levi-curvature assumption?
\end{question}

The purpose of this work is to answer Question~\ref{qII}.

\subsection{Our main results}\label{s-1.1}

Let us now formulate our main results. We refer to section~\ref{s-1.2} for some standard notations and terminology used here. Let $(X,T^{1,0}X)$ be a compact connected CR manifold of dimension $2n-1$, $n\geqslant2$. Let $L$ be a rigid CR line bundle over $X$. For every $u\in\Omega^{0,q}(X,L^k)$, we can define $Tu\in\Omega^{0,q}(X,L^k)$ and we have
\begin{equation}\label{e-gue140903IV}
T\ddbar_b=\ddbar_bT\ \ \mbox{on $\Omega^{0,q}(X,L^k)$},
\end{equation}
where $\ddbar_b:\Omega^{0,q}(X,L^k)\To\Omega^{0,q+1}(X,L^k)$ denotes the tangential Cauchy-Riemann operator.
For every $m\in\mathbb Z$, put
\begin{equation}\label{e-gue140903V}
\Omega^{0,q}_m(X,L^k):=\set{u\in\Omega^{0,q}(X,L^k);\, Tu=imu}.
\end{equation}
From \eqref{e-gue140903IV}, we have the $\ddbar_b$-complex for every $m\in\mathbb Z$:
\begin{equation}\label{e-gue140903VI}
\ddbar_b:\cdots\To\Omega^{0,q-1}_m(X,L^k)\To\Omega^{0,q}_m(X,L^k)\To\Omega^{0,q+1}_m(X,L^k)\To\cdots.
\end{equation}
For every $m\in\mathbb Z$, the $m$-th Fourier component of  $\ddbar_b$ cohomology is given by
\begin{equation}\label{e-gue140903VIIa}
H^{q}_{b,m}(X,L^k):=\frac{{\rm Ker\,}\ddbar_{b}:\Omega^{0,q}_m(X,L^k)\To\Omega^{0,q+1}_m(X,L^k)}{\operatorname{Im}\ddbar_{b}:\Omega^{0,q-1}_m(X,L^k)\To\Omega^{0,q}_m(X,L^k)}.
\end{equation}
The starting point of this paper is that without any Levi curvature assumption, for every $m\in\mathbb Z$ and every
$q=0,1,2,\ldots,n-1$, we have
\begin{equation}\label{e-gue140903VII}
{\rm dim\,}H^{q}_{b,m}(X,L^k)<\infty.
\end{equation}
Fix $\lambda\geq0$ and set $H^{q}_{b,\leq\lambda}(X,L^k):=\bigoplus\limits_{m\in\mathbb Z,\abs{m}\leq\lambda}H^q_{b,m}(X,L^k)$. In this work, we study the asymptotic behavior of the space $H^q_{b,\leq k\delta}(X,L^k)$ and its partial Szeg\H{o} kernel. Our main results are the following

\begin{theorem}[weak Morse inequalities]\label{t-mainI}
For $k$ large and for every $q=0,1,2,\ldots,n-1$, we have
\begin{equation}\label{jjjj}
{\rm dim\,}H^q_{b,\leq k\delta}(X,L^k)\leq(2\pi)^{-n}\frac{(-1)^{q}}{(n-1)!}k^n\int_X\int_{\mathbb R_{x,q}\bigcap[-\delta,\delta]}(i\mathcal R^L_x+i2s\mathcal L_x)^{n-1}\wedge(-\omega_0(x))ds+o(k^n),
\end{equation}
where $\mathcal R^L_x$ denotes the curvature of $L$, $\mathcal{L}_x$ denotes the Levi form of $X$, $\omega_0$ is the unique global non-vanishing real one form determined by $\langle\,\omega_0\,,\,U\,\rangle=0$, $\forall U\in T^{1,0}X\oplus T^{0,1}X$ and $\langle\,\omega_0\,,\,T\,\rangle=-1$ and
\begin{equation}\label{e-I}
\mathbb R_{x,q}:=\{s\in\Real: \mbox{$\mathcal R^L_x+2s\mathcal L_x$ has exactly $q$ negative and $n-1-q$ positive eigenvalues}\}.
\end{equation}
\end{theorem}
Although the eigenvalues of the Hermitian quadratic form $\mathcal R_x^L+2s\mathcal L_x, s\in\mathbb R$ are calculated with respect to the rigid Hermitian metric $\langle\cdot|\cdot\rangle$, the sign does not depend on the metric.
Note that $\mathcal R^L_x, \mathcal{L}_x\in T^{*1,0}_xX\wedge T^{*0,1}_xX$ (see  Definition~\ref{d-1.2}). Hence, $(\mathcal R^L_x+2s\mathcal L_x)^{n-1}\wedge(-\omega_0(x))$ is a global $2n-1$ form on $X$. Any Hermitian fiber metric $h^L$ on $L$ induces a curvature $\mathcal{R}^L$. It is easy to see that the integral in (\ref{jjjj}) does not depend on the choice of Hermitian fiber metric of $L$.

\begin{theorem}[strong Morse inequalities]\label{t-mainII}
For $k$ large and for every $q=0,1,2,\ldots,n-2$, we have
\begin{equation}\label{jjjjII}
\begin{split}
&\sum^q_{j=0}(-1)^{q-j}{\rm dim\,}H^j_{b,\leq k\delta}(X,L^k)\\
&\leq(2\pi)^{-n}\frac{k^n}{(n-1)!}(-1)^{q}\sum^q_{j=0}\int_X\int_{\mathbb R_{x,j}\bigcap[-\delta,\delta]}(i\mathcal R^L_x+i2s\mathcal L_x)^{n-1}\wedge(-\omega_0(x))ds+o(k^n),
\end{split}
\end{equation}
and when $q=n-1$, we have asymptotic Riemann-Roch-Hirzebruch theorem
\begin{equation}\label{jjjjIII}
\begin{split}
&\sum^{n-1}_{j=0}(-1)^{j}{\rm dim\,}H^j_{b,\leq k\delta}(X,L^k)\\
&=(2\pi)^{-n}\frac{k^n}{(n-1)!}\sum^{n-1}_{j=0}\int_X\int_{\mathbb R_{x,j}\bigcap[-\delta,\delta]}(i\mathcal R^L_x+i2s\mathcal L_x)^{n-1}\wedge(-\omega_0(x))ds+o(k^n).
\end{split}
\end{equation}
\end{theorem}

Demailly \cite{D85, D91} proved remarkable asymptotic Morse inequalities for the $\overline\partial$ complex constructed over the line bundle $L^k$ on compact complex manifold as $k\rightarrow\infty$, where $L$ is a holomorphic Hermitian line bundle. He solved with their help a generalized version of the Grauert-Riemenschneider. The original version of the conjecture had been solved previously by Siu \cite{S84, S85}. Shortly after, Bismut\cite{Bi87} gave a heat equation proof of Demailly's inequalities which involves probability theory.
\begin{definition}
We say that $(L, h^L)$ is a  positive rigid CR line bundle over $X$ if for any point $p\in X$, $\mathcal R_p^L$ is a positive Hermitian quadratic over $T_p^{1,0}X.$
\end{definition}
Assume that $\mathcal R^L$ is positive. The point of this paper is that if $\delta>0$ is small enough then $\mathbb R_{x,j}\cap [\delta, \delta]=\emptyset$, $\forall x\in X$ and for every $j=1,2,\ldots,n-1$. From this observation, \eqref{jjjj} and \eqref{jjjjIII}, we conclude that
\[{\rm dim\,}H^0_{b,\leq k\delta}(X,L^k)=(2\pi)^{-n}\frac{1}{(n-1)!}k^n\int_X\int_{\mathbb R_{x,0}\bigcap[-\delta,\delta]}(i\mathcal R^L_x+i2s\mathcal L_x)^{n-1}\wedge(-\omega_0(x))ds+o(k^n).\]
Hence, ${\rm dim\,}H^0_{b,\leq k\delta}(X,L^k)\approx k^n$. We conclude that

\begin{theorem}\label{llll}
If $L$ is a positive rigid CR line bundle, then $L$ is big, that is $\dim H_b^0(X, L^k)\gtrsim k^n$ when $k\gg1$.
\end{theorem}

We notice that from Theorem~\ref{t-mainI} and Theorem~\ref{t-mainII} and some simple argument, we can easily deduce Demailly's weak and strong Morse inequalities  (see the proof of Corollary~\ref{kkkk}).

\begin{definition}\label{d-I}
We say that condition $X(q)$ holds on $X$ if there is a $\delta>0$ such that $\mathbb R_{x,q}\bigcap[-\delta,\delta]=\emptyset$, $\forall ~x\in X$.
\end{definition}

In this work, we generalize Grauert-Riemenschneider criterion to CR manifolds with $S^1$ action and to general $(0,q)$-forms.

\begin{theorem}[Grauert-Riemenschneider criterion]\label{t-mainIII}
Given $q\in\{0,1,\ldots,n-1\}$, assume that $X(q-1)$ and $X(q+1)$ hold on $X$. Then, for some $\delta>0$,
\begin{equation}\label{l}
{\rm dim\,}H^q_{b,\leq k\delta}(X,L^k)=(2\pi)^{-n}\frac{(-1)^{q}}{(n-1)!}k^n\int_X\int_{\mathbb R_{x,q}\bigcap[-\delta,\delta]}(i\mathcal R^L_x+i2s\mathcal L_x)^{n-1}\wedge(-\omega_0(x))ds+o(k^n).
\end{equation}
\end{theorem}

\begin{definition}
We say that $L$ is a semi-positive rigid CR line bundle over $X$ if there exists a constant $\delta>0$ such that $\mathcal R_x^L+2s\mathcal L_x$ is a semi-positive Hermitian quadratic over $T_x^{1,0}X$ for any $x\in X,$ $|s|<\delta.$
\end{definition}

When $L$ is semi-positive, it is easy to see that condition $X(1)$ holds on $X$. From this observation and Theorem~\ref{t-mainIII}, we obtain the  Grauert-Riemenschneider criterion in the CR setting.

\begin{theorem}
If $L$ is a semi-positive rigid  CR line bundle and positive at a point, then $L$ is big.
\end{theorem}

\subsection{Set up and terminology}\label{s-1.2}
Let $(X, T^{1,0}X)$ be a compact connected CR manifold of dimension $2n-1, n\geq 2$, where $T^{1,0}X$ is a CR structure of $X$. That is, $T^{1,0}X$ is a subbundle of rank $n-1$ of the complexified tangent bundle $\mathbb{C}TX$, satisfying $T^{1,0}X\cap T^{0,1}X=\{0\}$, where $T^{0,1}X=\overline{T^{1,0}X}$, and $[\mathcal V,\mathcal V]\subset\mathcal V$, where $\mathcal V=C^\infty(X, T^{1,0}X)$. We assume that $X$ admits a $S^1$ action: $S^1\times X\rightarrow X$. We use $e^{i\theta}$ to denote the $S^1$ action. For $x\in X$, we say that the period of $x$ is $\frac{2\pi}{\ell}$, $\ell\in\mathbb N$, if $e^{i\theta}\circ x\neq x$, for every $0<\theta<\frac{2\pi}{\ell}$ and $e^{i\frac{2\pi}{\ell}}\circ x=x$. For each $\ell\in\mathbb N$, put
\begin{equation}\label{e-gue150802bm}
X_\ell=\set{x\in X;\, \mbox{the period of $x$ is $\frac{2\pi}{\ell}$}}
\end{equation}
and let $p=\min\set{\ell\in\mathbb N;\, X_\ell\neq\emptyset}$. It is well-known that if $X$ is connected, then $X_p$ is an open and dense subset of $X$ (see Duistermaat-Heckman~\cite{Du82} and Appendix in~\cite{HL15}) and the Lebesgue measure $m(X\setminus X_p)=0$. For simplicity, in this work, we assume that $p=1$ and we denote $X_{{\rm reg\,}}:=X_{1}$.

Let $T\in C^\infty(X, TX)$ be the global real vector field induced by the $S^1$ action given as follows
\begin{equation}\label{definition of Tu}
(Tu)(x)=\frac{\partial}{\partial\theta}\left(u(e^{i\theta}x)\right)\Big|_{\theta=0}, u\in C^\infty(X).
\end{equation}
\begin{definition}
We say that the $S^1$ action $e^{i\theta}, 0\leq\theta<2\pi$, is CR if
$$[T, C^\infty(X, T^{1,0}X)]\subset C^\infty(X, T^{1,0}X).$$
Furthermore, we say that the $S^1$ action is transversal if for each $x\in X$,
$$T(x)\oplus T_x^{1,0}(X)\oplus T_x^{0,1}X=\mathbb CT_xX.$$
\end{definition}
We assume throughout that $(X, T^{1,0}X)$ is a CR manifold with a transversal CR $S^1$ action $e^{i\theta}, 0\leq\theta<2\pi$ and we let $T$ be the global vector field induced by the $S^1$ action. Let $\omega_0\in C^\infty(X,T^*X)$ be the global real one form determined by $\langle\,\omega_0\,,\,U\,\rangle=0$, for every $U\in T^{1,0}X\oplus T^{0,1}X$ and $\langle\,\omega_0\,,\,T\,\rangle=-1$.
\begin{definition}\label{d-1.2}
For $x\in X$, the Levi-form $\mathcal L_x$ is the Hermitian quadratic form on $T_x^{1,0}X$ defined as follows. For any $U, V\in T_x^{1,0}X$, pick $\mathcal U, \mathcal V\in C^\infty(X, T^{1,0}X)$ such that $\mathcal U(x)=U, \mathcal V(x)=V$. Set
\begin{equation}
\mathcal L_x(U, \overline V)=\frac{1}{2i}\langle[\mathcal U, \overline{\mathcal V}](x), \omega_0(x)\rangle
\end{equation}
where $[\,,\,]$ denotes the Lie bracket. Note that $\mathcal L_x$ does not depend on the choice of $\mathcal U$ and $\mathcal V$.
\end{definition}
Denote by $T^{\ast 1,0}X$ and $T^{\ast0,1}X$ the dual bundles of
$T^{1,0}X$ and $T^{0,1}X$, respectively. Define the vector bundle of $(0,q)$-forms by
$T^{\ast 0,q}X:=\Lambda^qT^{\ast0,1}X$. Let $D\subset X$ be an open subset. Let $\Omega^{0,q}(D)$
denote the space of smooth sections of $T^{\ast0,q}X$ over $D$ and let $\Omega_0^{0,q}(D)$
be the subspace of $\Omega^{0,q}(D)$ whose elements have compact support in $D$. Similarly, if $E$ is a vector bundle, then we let $\Omega^{0,q}(D, E)$ denote the space of smooth sections of $T^{\ast0,q}X\otimes E$ over $D$ and let $\Omega_0^{0,q}(D, E)$ be the subspace of $\Omega^{0,q}(D, E)$ whose elements have compact support in $D$.

Fix $\theta_0\in [0, 2\pi)$. Let
$$d e^{i\theta_0}: \mathbb CT_x X\rightarrow \mathbb CT_{e^{i\theta_0}x}X$$
denote the differential map of $e^{i\theta_0}: X\rightarrow X$. By the property of transversal CR $S^1$ action, we can check that
\begin{equation}\label{a}
\begin{split}
de^{i\theta_0}:T_x^{1,0}X\rightarrow T^{1,0}_{e^{i\theta_0}x}X,\\
de^{i\theta_0}:T_x^{0,1}X\rightarrow T^{0,1}_{e^{i\theta_0}x}X,\\
de^{i\theta_0}(T(x))=T(e^{i\theta_0}x).
\end{split}
\end{equation}
Let $(de^{i\theta_0})^\ast: \Lambda^q(\mathbb CT^\ast X)\rightarrow\Lambda^q(\mathbb CT^\ast X)$ be the pull back of $de^{i\theta_0}, q=0,1\cdots, n-1$. From (\ref{a}), we can check that for  $q=0, 1,\cdots, n-1$,
\begin{equation}
(de^{i\theta_0})^\ast: T^{\ast0,q}_{e^{i\theta_0}x}X\rightarrow T_x^{\ast0,q}X.
\end{equation}
Let $u\in\Omega^{0,q}(X)$ and define $Tu$ as follows. For any $X_1,\cdots, X_q\in T_x^{1,0}X$,
\begin{equation}\label{b}
Tu(X_1,\cdots, X_q):=\frac{\partial}{\partial\theta}\left((de^{i\theta})^\ast u(X_1,\cdots, X_q)\right)\Big|_{\theta=0}.
\end{equation}
From the definition of $Tu$ it is easy to check that
$Tu=L_Tu$ for $u\in\Omega^{0, q}(X)$, where $L_Tu$ is the Lie derivative
of $u$ along the direction $T$. It is straightforward to see that (see also the discussion after Theorem~\ref{j})
\begin{equation}
T\overline\partial_b=\overline\partial_bT~\text{on}~\Omega^{0,q}(X).
\end{equation}
\begin{definition}
Let $D\subset X$ be an open set. We say that a function $u\in C^\infty(D)$ is rigid if $Tu=0$. We say that a function $u\in C^\infty(X)$ is Cauchy-Riemann (CR for short)
if $\ddbar_bu=0$. We say that $u\in C^\infty(X)$ is rigid CR if  $\ddbar_bu=0$ and $Tu=0$.
\end{definition}

\begin{definition}
Let $E$ be a complex vector bundle over $X$. We say that $E$ is rigid (resp.\ CR, resp.\ rigid CR) if there exists
an open cover $(U_j)_j$ of $X$ and trivializing frames $\set{f^1_j,f^2_j,\dots,f^r_j}$ on $U_j$,
such that the corresponding transition matrices are rigid (resp.\ CR, resp.\ rigid CR). The frames $\set{f^1_j,f^2_j,\dots,f^r_j}$ are called rigid (resp.\ CR, resp.\ rigid CR) frames.
\end{definition}

\begin{example}
Let $X$ be a compact CR manifold  with a locally free transversal CR $S^1$ action. Let $\{Z_j\}_j$ be a trivializing frame of $T^{1, 0}X$ defined in (\ref{e-can}). It is easy to check that the transition functions of such frames are rigid CR  and thus $T^{1, 0}X$ is a rigid CR vector bundle. Moreover, $\det{T^{1, 0}X}$ the determinant bundle of $T^{1, 0}X$ is a rigid CR line bundle.
\end{example}

\begin{example}
Let $(L, h)\overset{\pi}{\rightarrow}M$ be a Hermitian line bundle over a complex manifold $M$. Consider the circle bundle $X=\{v\in L: h(v)=1\}$ over $M$. Then $X$ is a compact CR manifold with a globally free transversal CR $S^1$ action. Let $E$ be a holomorphic vector bundle over $M$. Then the restriction of the pull back $\pi^\ast E|_X$ on $X$ is a rigid CR vector bundle over $X$.
\end{example}

From now on, let $L$ be a rigid CR line bundle over $X$.
We fix an open covering $(U_j)_j$ and a family $(s_j)_j$ of rigid CR frames $s_j$ on
$U_j$.
Let $L^k$ be the $k$-th tensor power of $L$.
Then $(s_j^{\otimes k})_j$ are rigid CR frames for $L^k$.
Let $s$ be a  rigid CR frame of $L$ on an open subset $D\subset X$ and locally for any $u\in\Omega^{0,q}(X, L)$, write $u=\tilde u\otimes s,$ $\tilde u\in\Omega^{0,q}(D)$, we define $Tu=T\tilde u\otimes s.$
Since the transition functions are rigid CR, $Tu$ is well defined. Moreover, we have
\begin{equation}\label{c}
T\overline\partial_b=\overline\partial_bT~\text{on}~\Omega^{0,q}(X, L).
\end{equation}

Fix a Hermitian fiber metric $h^L$ on $L$. If $s$ is a local rigid CR frame of $L$ on an open subset $D\subset X$, then the local weight of $h^L$ with respect to $s$ is the function $\Phi\in C^\infty(D, \mathbb R)$ for which
\begin{equation}\label{2016-06-27e1}
|s(x)|^2_{h^L}=e^{-\Phi(x)}, x\in D.
\end{equation}

\begin{definition}
Let $L$ be a rigid CR line bundle and let $h^L$ be a Hermitian metric on $L$.
The curvature of $(L,h^L)$ is the the Hermitian quadratic form $R^L=R^{(L,h)}$ on $T^{1,0}X$
defined by
\begin{equation}\label{e-gue150808w}
R_p^L(U,\overline V)=\frac12\,\big\langle d(\overline\partial_b\Phi-\partial_b\Phi)(p),
U\wedge\overline V\,\big\rangle,\:\: U, V\in T_p^{1,0}X,\:\: p\in D.
\end{equation}
\end{definition}
Due to \cite[Proposition 4.2]{HM12}, $R^L$ is a well-defined global Hermitian form,
since the transition functions between different rigid CR frames are annihilated by $T$.

\subsection{Hermitian CR geometry}
Fix a smooth Hermitian metric $\langle\cdot|\cdot\rangle$ on $\mathbb{C}TX$ so that
$T^{1,0}X$ is orthogonal to $T^{0,1}X$, $T$ is orthogonal to $T^{1,0}X\oplus T^{0,1}X$ and
$\langle T|T\rangle=1$. The Hermitian metric $\langle\cdot|\cdot\rangle$ on
$\mathbb CTX$ induces by duality a Hermitian metric on $\mathbb CT^\ast X$ and also on the bundles of $(0,q)$-forms $T^{\ast 0,q}X, q=0,1\cdots,n-1.$ We shall also denote all these induced
metrics by $\langle\cdot|\cdot\rangle$. For every $v\in T^{\ast0,q}X$, we write
$|v|^2:=\langle v|v\rangle$. We have the pointwise orthogonal decompositions
\begin{equation}
\begin{split}
&\mathbb CT^{\ast}X=T^{\ast1,0}X\oplus T^{\ast0,1}X\oplus\{\lambda\omega_0:\lambda\in\mathbb C\},\\
&\mathbb CTX=T^{1,0}X\oplus T^{0,1}X\oplus\{\lambda T:\lambda\in\mathbb C\}.
\end{split}
\end{equation}

\begin{definition}
Let $D$ be an open set and let $V\in C^\infty(D, \mathbb CTX)$ be a vector on $D$. We say that $V$ is rigid if
\begin{equation}
de^{i\theta}(V(x))=V(e^{i\theta}x)
\end{equation}
for any $x, \theta\in[0,2\pi)$ satisfying $x\in D, e^{i\theta}x\in D.$
\end{definition}
\begin{definition}
Let $\langle\cdot|\cdot\rangle$ be a Hermitian metric on $\mathbb CTX$. We say that $\langle\cdot|\cdot\rangle$ is rigid if for rigid vector fields $V, W$ on $D$, where $D$ is any on open set, we have
\begin{equation}
\langle V(x)|W(x)\rangle=\langle (de^{i\theta}V)(e^{i\theta}x)|(de^{i\theta}W)(e^{i\theta}x)\rangle, \forall x\in D, \theta\in[0,2\pi).
\end{equation}
\end{definition}
From theorem 9.2 in \cite{H14}, there is always a rigid Hermitian metric $\langle\cdot|\cdot\rangle$ on $\mathbb CTX$ such that $T^{1,0}X\bot T^{0,1}X, T\bot(T^{1,0}X\oplus T^{0,1}X), \langle T|T\rangle=1$ and $\langle u|v\rangle$ is real if $u, v$ are real tangent  vectors. Until further notice, we fix a rigid Hermitian metric $\langle\cdot|\cdot\rangle$ on $\mathbb CTX$ such that $T^{1,0}X\bot T^{0,1}X, T\bot(T^{1,0}X\oplus T^{0,1}X)$ and $ \langle T|T\rangle=1$.

\begin{definition}
Let $L$ be a rigid CR line bundle. A Hermitian fiber metric $h^L$ on $L$
is said to be rigid if
$T\Phi =0$ for local weight $\Phi$ with respect to any rigid CR frame.
\end{definition}
The definition does not depend on the choice of rigid CR frame.
\begin{lemma}\label{l-rigid}
There is a rigid Hermitian fiber metric on $L$. Moreover, for any Hermitian metric $\Td h^L$ on $L$, there is a rigid Hermitian metric $h^L$ of $L$ such that $\Td{\mathcal{R}}^L=R^L$ on $X$, where $\Td{\mathcal{R}}^L$ and $R^L$ denote the curvatures induced by $\Td h^L$ and $h^L$ respectively.
\end{lemma}
We will prove Lemma~\ref{l-rigid} in the end of section~\ref{q}.
Until furthermore, we assume that $h^L$ is a rigid Hermitian fiber metric on $L$. For $k>0,$ $ k\in\mathbb Z$, we shall consider $(L^k, h^{L^k})$. For $m\in\mathbb Z$, put
\begin{equation}
\Omega_m^{0,q}(X,L^k):=\{u\in\Omega^{0,q}(X, L^k): Tu=imu\}.
\end{equation}
Let $(\,\cdot\,|\,\cdot\,)_{h^{L^k}}$ be the $L^2$ inner product on $\Omega^{0,q}(X,L^k)$ induced by $h^{L^k}$, $\langle\,\cdot\,|\,\cdot\,\rangle$ and let $\norm{\cdot}_{h^{L^k}}$ denote the corresponding norm. Let $s$ be a local rigid CR frame of $L$ on an open set $D\subset X$. For $u=\Td u\otimes s^k, v=\Td v\otimes s^k\in\Omega^{0,q}_0(D,L^k)$, we have
\begin{equation}
(u|v)_{h^{L^k}}=\int_X\langle\tilde u|\tilde v\rangle e^{-k\Phi(x)}dv_X,
\end{equation}
where $dv_X$ is the volume form on $X$ induced by the rigid Hermitian metric $\langle\,\cdot\,|\,\cdot\,\rangle$.
Let $L^2_{(0,q),m}(X, L^k)$ be the completion of $\Omega_m^{0,q}(X,L^k)$ with respect to $(\cdot|\cdot)_{h^{L^k}}$. For $m\in\mathbb Z$, let
\begin{equation}
Q^q_{m,k}: L^2_{(0,q)}(X, L^k)\rightarrow L^2_{(0,q),m}(X, L^k)
\end{equation}
be the orthogonal projection with respect to $(\cdot|\cdot)_{h^{L^k}}$. Fix $\delta>0$, let $F_{\delta, k}:L^2_{(0,q)}(X, L^k)\rightarrow L^2_{(0,q)}(X, L^k)$ be the continuous map given by
\begin{equation}
F_{\delta,k}(u):=\sum\limits_{|m|\leq k\delta}Q^q_{m,k}u.
\end{equation}
Let $\overline\partial_{b,k}^\ast: \Omega^{0,q+1}(X, L^k)\rightarrow\Omega^{0,q}(X, L^k)$ be the formal adjoint of $\overline\partial_b$ with respect to $(\cdot|\cdot)_{h^{L^k}}$. Since $\langle\cdot|\cdot\rangle$ and $h^{L^k}$ are rigid, we can check that
\begin{equation}\label{d}
T\overline\partial_{b,k}^\ast=\overline\partial_{b,k}^\ast T~\text{on}~\Omega^{0,q}(X, L^k), q=0,1,\cdots,n-1,
\end{equation}
and
\begin{equation}\label{e}
\overline\partial_{b,k}^\ast:\Omega_m^{0,q+1}(X, L^k)\rightarrow\Omega_m^{0,q}(X, L^k),\forall m\in\mathbb Z.
\end{equation}
Put
$$\Box_{b,k}^{(q)}:=\overline\partial_b\overline\partial_{b,k}^\ast
+\overline\partial_{b,k}^\ast\overline\partial_b:\Omega^{0,q}(X, L^k)\rightarrow\Omega^{0,q}(X, L^k).$$
From (\ref{c}), (\ref{d})and (\ref{e}) we have
\begin{equation}
T\Box_{b,k}^{(q)}=\Box^{(q)}_{b,k}T~\text{on}~\Omega^{0,q}(X, L^k), q=0,1,\cdots,n-1,
\end{equation}
and
\begin{equation}
\Box_{b,k}^{(q)}:\Omega_m^{0,q}(X, L^k)\rightarrow\Omega_m^{0,q}(X, L^k),\forall m\in\mathbb Z.
\end{equation}
We will write $\Box^{(q)}_{b,k,m}$ to denote the restriction of $\Box^{(q)}_{b,k}$ on the space $\Omega_m^{0,q}(X, L^k)$. For every $m\in\mathbb Z$, we extend $\Box^{(q)}_{b,k,m}$ to $L^2_{(0,q),m}(X, L^k)$ in the sense of distribution by
\begin{equation}
\Box^{(q)}_{b,k,m}:{\rm Dom}(\Box^{(q)}_{b,k,m})\subset L^2_{(0,q),m}(X, L^k)\rightarrow L^2_{(0,q),m}(X, L^k),
\end{equation}
where ${\rm Dom}(\Box^{(q)}_{b,k,m})=\{u\in L^2_{(0,q),m}(X, L^k): \Box^{(q)}_{b,k,m}u\in L^2_{(0,q),m}(X, L^k)\}$.
The following follows from Kohn's $L^2$ estimate (see theorem 8.4.2 in \cite{CS01}).
\begin{theorem}\label{f}
For every $s\in\mathbb{N}_0=\mathbb N\cup\{0\}$, there exists a constant $C_{s,k}>0$ such that
\begin{equation}
\|u\|_{s+1}\leq C_{s,k}\left(\|\Box^{(q)}_{b,k}u\|_{s}+\|Tu\|_s+\|u\|_s\right),\forall u\in\Omega^{0,q}(X, L^k)
\end{equation}
where $\|\cdot\|_s$ denotes the sobolev norm of order $s$ on $X$.
\end{theorem}
From Theorem \ref{f}, we deduce that
\begin{theorem}\label{g}
Fix $m\in\mathbb Z$, for every $s\in\mathbb N_0$, there is a constant $C_{s,k,m}>0$ such that
\begin{equation}
\|u\|_{s+1}\leq C_{s,k,m}\left(\|\Box^{(q)}_{b,k,m}u\|_s+\|u\|_s\right), \forall u\in\Omega^{0,q}_m(X, L^k).
\end{equation}
\end{theorem}
From Theorem~\ref{g} and some standard argument in functional analysis, we deduce the following Hodge theory for $\Box^{(q)}_{b,k,m}$.
\begin{theorem}\label{gI}
Fix $m\in\mathbb Z$. $\Box^{(q)}_{b,k,m}:\mathrm{Dom}(\Box^{(q)}_{b,k,m})\subset L^2_{(0,q),m}(X, L^k)\rightarrow L^2_{(0,q),m}(X, L^k)$ is a self-adjoint operator. The spectrum of $\Box^{(q)}_{b,k,m}$ denoted by
$\mathrm {Spec}(\Box^{(q)}_{b,k,m})$ is a discrete subset of $[0,\infty)$. For every $\lambda\in\mathrm{Spec}(\Box^{(q)}_{b,k,m})$ the $\lambda$-eigenspace
\begin{equation}\label{h}
\mathcal{H}^q_{b, m,\lambda}(X, L^k):=\left\{u\in\mathrm{Dom}\Box^{(q)}_{b,k,m}:\Box^{(q)}_{b,k,m}u=\lambda u\right\}
\end{equation}
is finite dimensional with $\mathcal{H}^q_{b, m,\lambda}(X, L^k)\subset\Omega^{0,q}_m(X, L^k)$ and for $\lambda=0$ we denote by $\mathcal H^q_{b, m}(X, L^k)$ the harmonic space $\mathcal H^q_{b, m,0}(X, L^k)$ for brevity and then we have the Dolbeault isomorphism
\begin{equation}\label{i}
\mathcal{H}^q_{b, m}(X, L^k)\cong H^q_{b,m}(X, L^k).
\end{equation}
\end{theorem}
From Theorem \ref{gI} and (\ref{i}), we deduce that $\dim H^q_{b,m}(X, L^k)<\infty, \forall m\in\mathbb Z.$

\subsection{Our strategy}
Denote by $\det(\mathcal R_x^L+2s\mathcal L_x)$ the product of all the eigenvalues of $\mathcal R_x^L+2s\mathcal L_x$ with respect to the given rigid Hermitian metric.
Since
\[\Omega^{0, q}_m(X, L^k)\bot\Omega_{m^\prime}^{0, q}(X, L^k),\]
when $m, m^\prime \in\mathbb Z$ and $m\neq m^\prime$, we write
\begin{equation}
\Omega^{0, q}_{\leq k\delta}(X, L^k):=\bigoplus\limits_{m\in\mathbb Z, |m|\leq k\delta}\Omega^{0, q}_m(X, L^k)
\end{equation}
and in particular,
\begin{equation}
\mathcal H^q_{b, \leq k\delta}(X, L^k):=\bigoplus\limits_{m\in\mathbb Z, |m|\leq k\delta}\mathcal H^q_{b, m}(X, L^k).
\end{equation}
Here $\delta$ is a small constant. Then we have the following Hodge theory
\begin{equation}
\mathrm{dim}\mathcal H^q_{b, \leq k\delta}(X, L^k)<\infty, \mathcal H^q_{b,\leq k\delta}(X, L^k)\subset\Omega_{\leq k\delta}^{0,q}(X, L^k), \mathcal H^q_{b,\leq k\delta}(X, L^k)\cong H^q_{b,\leq k\delta}(X, L^k).
\end{equation}
Let $f_j\in\Omega^{0,q}_{\leq k\delta}(X, L^k), j=1,\cdots,m_k$ be an orthonormal basis for the space
$\mathcal H^q_{b, \leq k\delta}(X, L^k).$ The partial Szeg\H{o} kernel function is defined by
\begin{equation}
\Pi^{q}_{\leq k\delta}(x):=\sum_{j=1}^{m_k}|f_j(x)|^2_{h^{L^k}}.
\end{equation}
It is easy to see that $\Pi^{q}_{\leq k\delta}(x)$ is independent of the choice of the orthonormal basis and
\begin{equation}
\dim\mathcal H^q_{b, \leq k\delta}(X, L^k)=\int_X\Pi^{q}_{\leq k\delta}(x)dv_X.
\end{equation}
The following is our first main technique result.
\begin{theorem}\label{m}
\begin{equation}
\sup\{k^{-n}\Pi^{q}_{\leq k\delta}(x):k>0, x\in X\}<\infty.
\end{equation}
Furthermore, we have
\begin{equation}
\limsup_{k\rightarrow\infty}k^{-n}\Pi^{q}_{\leq k\delta}(x)\leq(2\pi)^{-n}\int_{\mathbb R_{x,q}\bigcap[-\delta,\delta]}|\det(\mathcal R_x^L+2s\mathcal L_x)|ds
\end{equation}
for all $x\in X$.
\end{theorem}
From Theorem \ref{m} and by Fatou's lemma, we obtain Theorem~\ref{t-mainI}.
From Theorem~\ref{t-mainI} and some simple argument, we deduce

\begin{corollary}[Demailly's weak morse inequalities]\label{kkkk}
Let $M$ be a compact Hermitian manifold with $\rm{dim}_{\mathbb C}M=n-1$ and $(L, h^L)$ be a Hermitian line bundle over $M$. Then $\forall q=0,1,2,\ldots,n-1$,
\begin{equation}
\dim H_{\overline\partial}^{q}(M, L^k)\leq k^{n-1}(2\pi)^{-(n-1)}\int_{M(q)}|\det \mathcal R^L_x|dv_M(x)+o(k^{n-1}),
\end{equation}
where $H_{\overline\partial}^{q}(M, L^k)$ denotes the $q$-th $\ddbar$-cohomology group with values in $L^k$, $dv_M$ is the induced volume form on $M$, $\mathcal R^L_x, x\in M$ is the Ricci curvature of the Hermitian line bundle $(L, h^L)$ and  $M(q)$ is a subset of $M$ where $\mathcal R^L_x$ has exactly $q$ negative eigenvalues and $n-1-q$ positive eigenvalues.
\end{corollary}

\begin{proof}
Let $X=M\times S^1$. Then $X$ is a Levi-flat CR manifold of $\rm{dim}_{\mathbb R} X=2n-1$ with $S^1$ action $e^{i\theta}$ and the global induced vector field is $T=\frac{\partial}{\partial\theta}$. Let $\pi_1: X=M\times S^1\rightarrow M$ be the natural projection. Then $L_1:=\pi_1^\ast L$ is naturally a rigid CR line bundle over $X$. It is easy to see that
\begin{equation}\label{e-dm}
\dim H_{\overline\partial}^{q}(M, L^k)=\frac{1}{2k+1}\dim H^q_{b,\leq k}(X,  L_1^k).
\end{equation}
From \eqref{jjjj}, we have
\begin{equation}\label{e-dmI}
\begin{split}
&\dim H^q_{b,\leq k}(X,L_1^k)\\
&\leq(2\pi)^{-n}k^n\int_{X}\int_{\mathbb R_{x,q}\cap[-1,1]}|\det\mathcal R^{L_1}_x+2s\mathcal L_x|ds dv_X+o(k^n)\\
&=(2\pi)^{-n}k^n\int_{M\times S^1}\int_{\mathbb R_{x,q}\cap[-1,1]}|\det \mathcal R^{L_1}_x|ds dv_X+o(k^n)\\
&=2(2\pi)^{-n}k^n\int_{M(q)\times S^1}|\det\mathcal R^L_x|dv_Mdv_{S^1}+o(k^n)\\
&=2(2\pi)^{-(n-1)}k^n\int_{M(q)}|\det\mathcal R^L_x|dv_M+o(k^n).
\end{split}
\end{equation}
From \eqref{e-dmI} and \eqref{e-dm}, we get the conclusion of the Corollary \ref{kkkk}.
\end{proof}
For $\lambda\geq 0$ and $\lambda\in\mathbb R$, we define
\begin{equation}
\mathcal H^q_{b, \leq k\delta, \lambda}(X, L^k):=\left\{u\in\Omega^{0,q}_{\leq k\delta}(X, L^k):\Box^{(q)}_{b,k}u=\lambda u\right\}
\end{equation}
and
\begin{equation}
\begin{split}
\mathcal H^q_{b, \leq k\delta, \leq k\sigma}(X, L^k):=&\bigoplus\limits_{\lambda\leq k\sigma}\mathcal H^q_{b, \leq k\delta,\lambda}(X, L^k).
\end{split}
\end{equation}
Set $\Pi^{q}_{\leq k\delta, \leq k\sigma}(x)=\sum_{j=1}^{d_k}|g_j(x)|^2$, where $\{g_j(x)\}_{j=1}^{d_k}\subset\Omega^{0, q}_{\leq k\delta}(X, L^k)$ is any orthonormal basis of the space $\mathcal H^q_{b, \leq k\delta, \leq k\sigma}(X, L^k).$ Our second main technique result is the following
\begin{theorem}\label{o}
For any sequence $v_k>0$ with $v_k\rightarrow0$ as $k\rightarrow\infty$, there exists a constant $C_0^\prime$ independent of $k$, such that
\begin{equation}
k^{-n}\Pi^{q}_{\leq k\delta, \leq kv_k}(x)\leq C_0^\prime
\end{equation}
for all $x\in X$. Moreover, there is a sequence $\mu_k>0, \mu_k\rightarrow0$ as $k\rightarrow\infty$, such that for any sequence $v_k>0$ with $\lim\limits_{k\rightarrow\infty}\frac{\mu_k}{v_k}=0$, we have
\begin{equation}\label{p}
\lim\limits_{k\rightarrow\infty}k^{-n}\Pi^{q}_{\leq k\delta, \leq kv_k}(x)=(2\pi)^{-n}\int_{\mathbb R_{x,q}\bigcap[-\delta,\delta]}|\det(\mathcal R^L_x+2s\mathcal L_x)|ds
\end{equation}
for all $x\in X_{{\rm reg\,}}$.
\end{theorem}
Integrating (\ref{p}), we have
\begin{theorem}\label{xx}
There is a sequence $\mu_k>0, \mu_k\rightarrow0$ as $k\rightarrow\infty$, such that for any sequence $v_k>0$ with $\lim\limits_{k\rightarrow\infty}\frac{\mu_k}{v_k}=0$, we have
\begin{equation}
\dim\mathcal H^q_{b, \leq k\delta,\leq kv_k}(X, L^k)=(2\pi)^{-n}k^n\int_X\int_{\mathbb R_{x,q}\bigcap[-\delta,\delta]}|\det(\mathcal R^L_x+2s\mathcal L_x)|dsdv_X+o(k^n).
\end{equation}
\end{theorem}
\begin{proof}[Proof of Theorem~\ref{t-mainIII}]
Set $\mathcal H^q_{b, \leq k\delta, 0<\lambda\leq k\sigma}(X, L^k):=\bigoplus\limits_{ 0<\lambda\leq k\sigma}\mathcal H^q_{b, \leq k\delta,\lambda}(X, L^k).$ We define a map
\begin{equation}
\begin{split}
P: \mathcal H^q_{b, \leq k\delta, 0<\lambda\leq kv_k}(X, L^k)&\rightarrow\mathcal H^{q-1}_{b, \leq k\delta, 0<\lambda\leq kv_k}(X, L^k)\oplus\mathcal H^{q+1}_{b, \leq k\delta, 0<\lambda\leq kv_k}(X, L^{k})\\
u&\mapsto(\overline\partial_b^\ast u, \overline\partial_b u).
\end{split}
\end{equation}
Since map $P$ is injective, it follows that
\begin{equation}
{\rm dim}\mathcal H^q_{b, \leq k\delta, 0<\lambda\leq kv_k}(X, L^k)\leq{\rm dim}\mathcal H^{q-1}_{b, \leq k\delta, 0<\lambda\leq kv_k}(X, L^k)+{\rm dim}\mathcal H^{q+1}_{b, \leq k\delta, 0<\lambda\leq kv_k}(X, L^{k}).
\end{equation}
From Theorem~\ref{xx}, we have
\begin{equation}\label{mmmm}
{\rm dim}\mathcal H^{q-1}_{b, \leq k\delta, 0<\lambda\leq kv_k}(X, L^k)=o(k^n), {\rm dim}\mathcal H^{q+1}_{b, \leq k\delta, 0<\lambda\leq kv_k}(X, L^{k})=o(k^n).
\end{equation}
Since ${\rm dim}\mathcal H^q_{b, \leq k\delta, \leq kv_k}(X, L^k)={\rm dim}\mathcal H^q_{b, \leq k\delta, 0<\lambda\leq kv_k}(X, L^k)+{\rm dim}\mathcal H^q_{b, \leq k\delta}(X, L^k),$ combining Theorem~\ref{xx} and (\ref{mmmm}), we get the conclusion of the Theorem~\ref{t-mainIII}.
\end{proof}
From Theorem \ref{xx} and the linear algebraic argument from Demailly in \cite{D85}, \cite{D91} and \cite{M96}, we obtain Theorem~\ref{t-mainII}.
From Theorem~\ref{t-mainII}, we can repeat the proof of Corollary~\ref{kkkk} and deduce

\begin{corollary}[Demailly's strong Morse inequalities]
Let $M$ be a compact Hermitian manifold with $\rm{dim}_{\mathbb C}M=n-1$ and $(L, h^L)$ be a Hermitian line bundle on $M$. Then for any $0\leq q\leq n-1$, we have
\begin{equation}
\sum_{j=0}^q(-1)^{q-j}\dim H_{\overline\partial}^{j}(M, L^k)\leq k^{n-1}(2\pi)^{-(n-1)}\sum_{j=0}^q(-1)^{q-j}\int_{M(j)}|\det\mathcal R^L_x|dv_M+o(k^{n-1}).
\end{equation}
\end{corollary}

\subsection{Canonical local coordinates}\label{q}
In this work, we need the following result due to Baouendi-Rothschild-Treves, (see \cite{BRT85}).
\begin{theorem}\label{j}
For  $x_0\in X$, there exist local coordinates $(x_1,\cdots,x_{2n-1})=(z,\theta)=(z_1,\cdots,z_{n-1},\theta), z_j=x_{2j-1}+ix_{2j},j=1,\cdots,n-1, x_{2n-1}=\theta$, defined in some small neighborhood $D=\{(z, \theta): |z|<r, |\theta|<\varepsilon\}$ centered at $x_0$ such that on $D$
\begin{equation}\label{e-can}
\begin{split}
&T=\frac{\partial}{\partial\theta}\\
&Z_j=\frac{\partial}{\partial z_j}+i\frac{\partial\varphi}{\partial z_j}(z)\frac{\partial}{\partial\theta},j=1,\cdots,n-1,
\end{split}
\end{equation}
where $\{Z_j(x)\}$ form a basis of $T_x^{1,0}X$ for each $x\in D$, and $\varphi(z)\in C^\infty(D,\mathbb R)$ is independent of $\theta$.
\end{theorem}
The local coordinates defined in Theorem \ref{j} are called canonical local coordinates.
By using canonical local coordinates, we get another way to define $Tu, \forall u\in\Omega^{0,q}(X)$. Let $(z,\theta)$ be the canonical coordinates  defined on $D\subset\subset X$. It is clearly that
$$\{d\overline z_{j_1}\wedge\cdots\wedge d\overline z_{j_q}, 1\leq j_1<\cdots<j_q\leq n-1\}$$
is a basis for $T^{\ast0,q}_xX$, $\forall ~x\in D$. Let $u\in\Omega^{0,q}(X)$. On $D$, we write
$$u=\sum\limits_{j_1<\cdots<j_q}u_{j_1\cdots j_q}d\overline z_{j_1}\wedge\cdots\wedge d\overline z_{j_q}.$$
Then on $D$ we can check that
\begin{equation}\label{lI}
Tu=\sum\limits_{j_1<\cdots<j_q}(Tu_{j_1\cdots j_q})d\overline z_{j_1}\wedge\cdots\wedge d\overline z_{j_q}.
\end{equation}

\begin{remark}\label{k}
Since the Hermitian metric $\langle\cdot|\cdot\rangle$ on $\mathbb CTX$ is rigid, we can find orthonormal frame $\{e^j\}_{j=1}^{n-1}$ of $T^{\ast0,1}X$ on $D$ such that $e^j(x)=e^j(z), \forall x=(z, \theta)\in D,  j=1,\cdots,n-1$. Moreover, if we denote by $dv_X$ the volume form with respect to the rigid Hermitian metric on $\mathbb CTX$, then on $D$,
\begin{equation}
dv_X=m(z)dv(z)d\theta,
\end{equation}
where $m(z)\in C^\infty(D, \mathbb R)$ which does not depend on $\theta$ and $dv(z)=2^{n-1}dx_1\cdots dx_{2n-2}$.
\end{remark}
With respect to the orthonormal frame defined in Remark \ref{k}, write
\begin{equation}
u=\sum\nolimits_{|J|=q}^\prime u_Je^J, J=(j_1, \ldots, j_q), e^J=e^{j_1}\wedge\cdots\wedge e^{j_q},
\end{equation}
where the prime means the multi index in the summation is strictly increasing.
Then from (\ref{l}) and Remark \ref{k}, we can check that
\begin{equation}
Tu=\sum\nolimits_{|J|=q}^\prime(Tu_J)e^J.
\end{equation}
\begin{proof}[Proof of Lemma~\ref{l-rigid}]
Fix $p\in X$ and let $(z,\theta)$ be canonical coordinates defined in some neighbourhood of $p$ such that $(z(p),\theta(p))=(0,0)$ and \eqref{e-can} hold. Suppose that $(z,\theta)$ defined on $\{z\in\mathbb C^{n-1}:\, \abs{z}<\delta\}\times\{\theta\in\Real:\, \abs{\theta}<\delta\}$, for some $\delta>0$. For $z\in\mathbb C^{n-1}$, $\abs{z}<\delta$, $\theta\in\Real$, we identify $(z,\theta)$ with $e^{i\theta}\circ (z,0)\in X$. Thus, we may assume that $\theta$ is defined on $\Real$. Put
\[\begin{split}
A:=&\{\lambda\in[0,2\pi]:\, \mbox{There is a local trivializing section $s$ defined on}\\
&\quad\quad\quad\mbox{$\set{z\in\Complex^{n-1}:\,\abs{z}<\varepsilon}\times[0,\lambda+\varepsilon)$, for some $0<\varepsilon<\delta$}\}.\end{split}\]
It is clearly that $A$ is a non-empty open set in $[0,2\pi]$. We claim that $A$ is closed. Let $\lambda_0$ be a limit point of $A$. Consider the point $(0,\lambda_0)$. For some $\varepsilon_1>0$, $\varepsilon_1$ small, there is a local trivializing section $s_1$ defined on $\set{z\in\Complex^{n-1}:\,\abs{z}<\varepsilon_1}\times(\lambda_0-\varepsilon_1,\lambda_0+\varepsilon_1]$. Since $\lambda_0$ is a limit point of $A$, we can find a local trivializing section $\Td s$ defined on $\set{z\in\Complex^{n-1}:\,\abs{z}<\varepsilon_2}\times[0,\lambda_0-\frac{\varepsilon_1}{2})$, for some $\varepsilon_2>0$. Now, $\Td s=gs_1$ on
\[\set{z\in\Complex^{n-1}:\,\abs{z}<\varepsilon_0}\times(\lambda_0-\varepsilon_1,\lambda_0-\frac{\varepsilon_1}{2})\]
for some rigid CR function $g$, where $\varepsilon_0=\min\set{\varepsilon_1,\varepsilon_2}$. Since $g$ is independent of $\theta$, $g$ is well-defined on $\set{z\in\Complex^{n-1}:\,\abs{z}<\varepsilon_0}\times\Real$. Put $s=\Td s$ on $\set{z\in\Complex^{n-1}:\,\abs{z}<\varepsilon_0}\times[0,\lambda_0-\frac{\varepsilon_1}{2})$ and $s=gs_1$ on $\set{z\in\Complex^{n-1}:\,\abs{z}<\varepsilon_0}\times[\lambda_0-\frac{\varepsilon_1}{2},\lambda_0+\varepsilon_1)$. It is straightforward to check that $s$ is well-defined as a local trivializing section on $\set{z\in\Complex^{n-1}:\,\abs{z}<\varepsilon_0}\times[0,\lambda_0+\varepsilon_1)$. Thus, $\lambda_0\in A$ and hence $A=[0,2\pi]$.

From the discussion above, we see that we can find local trivializations $W_1,\ldots,W_N$ such that $X=\bigcup^N_{j=1}W_j$ and $\bigcup_{0\leq\theta\leq2\pi}e^{i\theta}W_t\subset W_t$, $t=1,\ldots,N$.  Take any
Hermitian fiber metric $\Td h^{L}$ on $L$ and let $\Td\Phi$ denotes the corresponding local weight. Let $h^{L}$ be the Hermitian fiber metric on $L$ locally given by
$\abs{s}^2_{h^{L}}=e^{-\Phi}$, where $\Phi(x)=\frac{1}{2\pi}\int^{2\pi}_{0}\Td\Phi(e^{i\theta}x)d\theta$. It is obviously that $h^{L}$
is well-defined and $T\Phi=0$. Moreover, it is easy to see that $\Td{\mathcal{R}}^L=\mathcal{R}^L$, where $\Td{\mathcal{R}}^L$ and $\mathcal{R}^L$ denote the curvatures of $L$ induced by $\Td h^L$ and $h^L$ respectively. The lemma follows.
\end{proof}

\section{The estimates of the partial Szeg\H{o} kernel function $\Pi_{\leq k\delta}^{q}$}

We first introduce some notations. For $x\in X$, we can choose an orthonormal frame $\{e^j\}_{j=1}^{n-1}$ of $T^{\ast0,1}X$ defined in section \ref{q} over a neighborhood $D$ of $x$. For $J=(j_1,\cdots, j_q)$ with $j_1<\cdots<j_q$, we define $e^J=e^{j_1}\wedge\cdots\wedge e^{j_q}$. Then $\{e^J:|J|=q, J~\text{strictly increasing}\}$ is an orthonormal frame for $T^{\ast0,q}X$ over $D$. For any $f\in\Omega^{0,q}(X, L^k)$, on $D$, we may write
\begin{equation}
f=\sum\nolimits_{|J|=q}^\prime f_Je^J, ~\text{with}~f_J=\langle f|e^J\rangle\in C^\infty(D, L^k).
\end{equation}
The extramal function $S^q_{\leq k\delta, J}(y)$ for $y\in D$ along the direction $e^J$ is defined by
\begin{equation}
S^q_{\leq k\delta, J}(y):=\sup\limits_{\alpha\in\mathcal H^q_{b, \leq k\delta}(X, L^k),  \|\alpha\|_{h^{L^k}}=1}|\alpha_J(y)|_{h^{L^k}}^2.
\end{equation}
We can repeat the proof of Lemma 2.1 in~\cite{HM12} and conclude that

\begin{lemma}\label{z}
For every local orthonormal frame $\{e^J:|J|=q,~\text{strict increasing}\}$ of $T^{\ast0,q}X$ over an open set $D$ we have for $y\in D$,
\begin{equation}
\Pi_{\leq k\delta}^q(y)=\sum\nolimits_{|J|=q}^\prime S^q_{\leq k\delta, J}(y).
\end{equation}
\end{lemma}

\subsection{The scaling technique}\label{kk}

Fix  $p\in X$. Let $U_1,\cdots, U_{n-1}$ be the dual frame of $\overline{e^1},\cdots, \overline{e^{n-1}}$ and for which the Levi-form is diagonal at $p$. Furthermore, let $s$ be a rigid CR frame of  $L$ on an open neighborhood of $p$ and $|s|^2_{h^L}=e^{-\Phi}$. We take canonical local coordinates $(z,\theta), z_j=x_{2j-1}+i x_{2j}, j=1,\cdots, n-1$ defined in Theorem \ref{j} such that $\omega_0(p)=-d\theta, (z(p), \theta(p))=0$,
\begin{equation}\label{e-gue150303}
\left\langle\frac{\partial}{\partial x_j}(p)\Big|\frac{\partial}{\partial x_t}(p)\right\rangle=2\delta_{jt}, \left\langle\frac{\partial}{\partial x_j}(p)\Big|\frac{\partial}{\partial \theta}(p)\right\rangle=0, \left\langle\frac{\partial}{\partial \theta}(p)\Big|\frac{\partial}{\partial \theta}(p)\right\rangle=1
\end{equation}
for $j,t=1,\ldots, 2n-2$, and with respect to the canonical coordinates $(z,\theta)$,
\begin{equation}\label{e-gue150303I}
U_j=\frac{\partial}{\partial z_j}+i\lambda_j\overline z_j\frac{\partial}{\partial\theta}+O(|(z,\theta)|^2)
,j=1,\ldots,n-1,
\end{equation}
where $\{\lambda_j\}_{j=1}^{n-1}$ are the eigenvalues of Levi-form at $p$ with respect to the given rigid Hermitian metric and $\frac{\partial}{\partial z_j}=\frac{1}{2}(\frac{\partial}{\partial{x_{2j-1}}}-i\frac{\partial}{\partial x_{2j}}), j=1,\ldots, n-1$. Moreover, by changing the local rigid CR frame of $L$ we assume the local weight
\begin{equation}\label{e-gue150303II}
\Phi=\sum\limits_{j,t=1}^{n-1}\mu_{j,t}\overline z_jz_t+O(|z|^3).
\end{equation}
In this section, we work with canonical local coordinates $(z,\theta)$ defined on an open neighbourhood $D$ of $p$ and we identify $D$ with some open set in $\mathbb R^{2n-1}$. Let $(\cdot|\cdot)_{k\Phi}$ be the weighted inner product on the space $\Omega^{0,q}_0(D)$ defined as follows:
\begin{equation}
(f|g)_{k\Phi}=\int_D\langle f|g\rangle e^{-k\Phi(z)}dv_X(x)
\end{equation}
where $f, g\in\Omega_0^{0,q}(D)$. We denote by $L^2_{(0,q)}(D,k\Phi)$  the completion of $\Omega_0^{0,q}(D)$ with respect to $(\cdot|\cdot)_{k\Phi}$. For $r>0$, let
$D_r=\{(z,\theta)\in\mathbb R^{2n-1}: |z|<r, |\theta|<r\}$. Here $\{z=(z_1,\ldots, z_{n-1})\in\mathbb
C^{n-1}:|z|<r\}$ means that $\{z\in\mathbb C^{n-1}:|z_j|<r, j=1,\cdots, n-1\}$. Let $F_k$ be the scaling map $F_k(z,\theta)=(\frac{z}{\sqrt k}, \frac{\theta}{k})$. From now on, we
assume $k$ is sufficiently large such that $F_k(D_{\log k})\Subset D$. We define the scaled bundle $F_k^\ast T^{\ast0, q} X$ on $D_{\log k}$ to be the bundle whose fiber at
$(z,\theta)\in D_{\log k}$ is
\begin{equation}
F_k^\ast T^{\ast0, q}X|_{(z,\theta)}=\left\{\sum\nolimits_{|J|=q}^\prime a_Je^J\left(\frac{z}{\sqrt k}, \frac{\theta}{k}\right): a_J\in\mathbb C, |J|=q, J~\text{strictly increasing}\right\}.
\end{equation}
We take the Hermitian metric $\langle\cdot|\cdot\rangle_{F_k^\ast}$ on $F_k^\ast T^{\ast 0, q} X$ so that at each point $(z,\theta)\in D_{\log k}$,
\begin{equation}
\left\{e^J\left(\frac{z}{\sqrt k}, \frac{\theta}{k}\right): |J|=q, J~\text{strictly increasing}\right\}
\end{equation}
is an orthonormal frame for $F_k^\ast T^{\ast0, q} X$ on $D_{\log k}$. Let $F_k^\ast\Omega^{0,q}(D_r)$ denote the space of smooth sections of $F_k^\ast T^{\ast0, q} X$ over $D_r$ and let $F_k^\ast\Omega^{0,q}_0(D_r)$ be the subspace of $F_k^\ast\Omega^{0,q}(D_r)$ whose elements have compact support in $D_r$.  Given $f\in\Omega^{0,q}(D_r)$.
We write $f=\sum\nolimits_{|J|=q}^\prime f_Je^J$. We define the scaled form $F_k^\ast f\in F_k^\ast\Omega^{0,q}(D_{\log k})$ by
\begin{equation}
F_k^\ast f=\sum\nolimits_{|J|=q}^\prime f_J\left(\frac{z}{\sqrt k}, \frac{\theta}{k}\right)e^J\left(\frac{z}{\sqrt k}, \frac{\theta}{k}\right).
\end{equation}
For brevity, we denote $F_k^\ast f$ by $f(\frac{z}{\sqrt k}, \frac{\theta}{k})$.
Let $P$ be a partial differential operator of order one on $F_k(D_{\log k})$ with $C^\infty$ coefficients. We write $P=a(z,\theta)\frac{\partial}{\partial\theta}+\sum\limits_{j=1}^{2n-2}a_j(z, \theta)\frac{\partial}{\partial x_j}.$ The scaled partial differential operator $P_{(k)}$ on $D_{\log k}$
is given by
\begin{equation}
P_{(k)}=\sqrt{k}F_k^\ast a\frac{\partial}{\partial\theta}+\sum_{j=1}^{2n-2}F_k^\ast a_j\frac{\partial}{\partial x_j}.
\end{equation}
Let $f\in C^\infty(F_k(D_{\log k}))$. We can check that
\begin{equation}
P_{(k)}(F_k^\ast f)=\frac{1}{\sqrt k}F_k^\ast(Pf).
\end{equation}
The scaled differential operator $\overline\partial_{b,(k)}: F_k^\ast\Omega^{0,q}(D_{\log k})\rightarrow F_k^\ast\Omega^{0,q+1}(D_{\log k})$ is given by
\begin{equation}\label{ee1}
\overline\partial_{b,(k)}=\sum_{j=1}^{n-1}e_j\left(\frac{z}{\sqrt k}, \frac{\theta}{k}\right)\wedge\overline U_{j,(k)}+\sum_{j=1}^{n-1}\frac{1}{\sqrt{k}}(\overline\partial_be_j)\left(\frac{z}{\sqrt k}, \frac{\theta}{k}\right)\wedge\left(e_j\left(\frac{z}{\sqrt k}, \frac{\theta}{k}\right)\wedge\right)^\ast,
\end{equation}
where $\left(e_j\left(\frac{z}{\sqrt k}, \frac{\theta}{k}\right)\wedge\right)^\ast:F^\ast_kT^{\ast0,q}X\To F^\ast_kT^{\ast0,q-1}X$ is the adjoint of $e_j\left(\frac{z}{\sqrt k}, \frac{\theta}{k}\right)\wedge$ with respect to the $\langle\cdot\big|\cdot\rangle_{F_k^\ast}$, $ j=1,\ldots, n-1$. That is, $$\left\langle e_j\left(\frac{z}{\sqrt{k}},\frac{\theta}{k}\right)\wedge u \Big|v \right\rangle_{F^\ast_k}=\left\langle u\Big|\left(e_j\left(\frac{z}{\sqrt{k}},\frac{\theta}{k}\right)
\wedge\right)^\ast v\,\right\rangle_{F^\ast_k}$$ for all $u\in F^\ast_kT^{\ast0,q-1}X$, $v\in F^\ast_kT^{\ast0,q}X$.
From (\ref{ee1}), $\overline\partial_{b,(k)}$ satisfies that
\begin{equation}\label{l1}
\overline\partial_{b,(k)}F_k^\ast f=\frac{1}{\sqrt{k}}F_k^\ast(\overline\partial_b f).
\end{equation}
Let $(\cdot|\cdot)_{kF_k^\ast\Phi}$ be the inner product on the space $F_k^\ast\Omega^{0,q}_0(D_{\log k})$ defined  as follows:
\begin{equation}
(f|g)_{kF_k^\ast\Phi}=\int_{D_{\log k}}\langle f|g\rangle_{F_k^\ast}e^{-kF_k^\ast\Phi}(F_k^\ast m)dv(z)d\theta.
\end{equation}
Let
\[\overline\partial^\ast_{b, (k)}: F_k^\ast\Omega^{0,q+1}(D_{\log k})\rightarrow F_k^\ast\Omega^{0,q}(D_{\log k})\]
be the formal adjoint of $\overline\partial_{b, (k)}$ with respect to $(\cdot|\cdot)_{kF_k^\ast\Phi}$. Then we also have
\begin{equation}\label{l2}
\overline\partial_{b,(k)}^\ast F_k^\ast f=\frac{1}{\sqrt{k}}F_k^\ast(\overline\partial_{b}f).
\end{equation}
We define now the scaled complex Laplacian $\Box^{(q)}_{b,(k)}:F_k^\ast\Omega^{0,q}(D_{\log k})\rightarrow F_k^\ast\Omega^{0,q}(D_{\log k})$ which is given by
\begin{equation}
\Box^{(q)}_{b,(k)}=\overline\partial_{b,(k)}^\ast\overline\partial_{b,(k)}
+\overline\partial_{b,(k)}\overline\partial_{b,(k)}^\ast.
\end{equation}
Then (\ref{l1}) and (\ref{l2}) imply that
\begin{equation}\label{s}
\Box^{(q)}_{b,(k)}F_k^\ast f=\frac{1}{k}F_k^\ast(\Box^{(q)}_{b,k}f).
\end{equation}
Similarly, as Proposition 2.3 in \cite{H08}, Proposition 2.3 in \cite{HM12},  we have
\begin{proposition}
\begin{equation}
\begin{split}
\Box^{(q)}_{b,(k)}=&\sum_{j=1}^{n-1}\left[\left(-\frac{\partial}{\partial z_j}-i\lambda_j\overline z_j\frac{\partial}{\partial\theta}+\sum_{t=1}^{n-1}\mu_{t,j}\overline z_t\right)\left(\frac{\partial}{\partial\overline z_j}-i\lambda_j z_j\frac{\partial}{\partial\theta}\right)\right]\\
&+\sum_{j,t=1}^{n-1}e_j\left(\frac{z}{\sqrt k},\frac{\theta}{k}\right)\wedge\left(e_t\left(\frac{z}{\sqrt k},\frac{\theta}{k}\right)\wedge\right)^\ast\left(\mu_{j,t}
-2i\lambda_j\delta_{j,t}\frac{\partial}{\partial\theta}\right)+\varepsilon_kP_k
\end{split}
\end{equation}
on $D_{\log k}$,
where $\varepsilon_k$ is a sequence tending to zero with $k\rightarrow\infty$, $P_k$ is a second order differential operator and all the derivatives of the coefficients of $P_k$ are uniformly bounded  in $k$ on $D_{\log k}$.
\end{proposition}
Let $U\subset D_{\log k}$ be an open set and let $W^s_{kF_k^\ast\Phi}(U, F_k^\ast T^{\ast0, q} X), s\in\mathbb N_0:=\mathbb N\cup\{0\}$ denote the sobolev space of order $s$ of sections of $F_k^\ast T^{\ast0,q} X$ over $U$ with respect to the weight $kF_k^\ast\Phi$. The sobolev norm on this space is given by
\begin{equation}
\|u\|^2_{kF_k^\ast\Phi,s,U}:=\sum\limits_{\alpha\in\mathbb N_0^{2n-1}, |\alpha|\leq s}\sum\nolimits_{|J|=q}^\prime
\int_{U}|\partial^\alpha_{x,\theta}u_J|^2e^{-kF_k^\ast\Phi}(F_k^\ast
m)dv(z)d\theta
\end{equation}
where $u=\sum\nolimits_{|J|=q}^\prime u_Je^J\left(\frac{z}{\sqrt k},\frac{\theta}{k}\right)\in W^s_{kF_k^\ast\Phi}(U, F_k^\ast T^{\ast0,q} X)$.

\begin{proposition}\label{t}
For every $r>0$ with $D_{2r}\Subset D_{\log k}$, there exists a constant $C_{r,s}>0$ independent of $k$ and the point $p$ such that for all $u\in F_k^\ast\Omega^{0, q}(D_{\log k})$, we have
\begin{equation}
\|u\|_{kF_k^\ast\Phi, s+1, D_r}^2\leq C_{r,s}\left(\|u\|^2_{kF_k^\ast\Phi, D_{2r}}+\|\Box^{(q)}_{b,(k)}u\|^2_{kF_k^\ast\Phi, s, D_{2r}}
+\left\|\left(\frac{\partial}{\partial\theta}\right)^{s+1}u\right\|^2_{kF_k^\ast\Phi,D_{2r}}\right).
\end{equation}
\end{proposition}
\begin{proof}
We can repeat the procedure of Kohn's $L^2$ estimate with minor change (see Theorem 8.4.2 in \cite{CS01}) and conclude that
\begin{equation}
\|u\|_{kF_k^\ast\Phi, s+1, D_r}^2\leq C_{r,s,k}\left(\|u\|^2_{kF_k^\ast\Phi, D_{2r}}+\|\Box^{(q)}_{b,(k)}u\|^2_{kF_k^\ast\Phi, s, D_{2r}}
+\left\|\left(\frac{\partial}{\partial\theta}\right)^{s+1}u\right\|^2_{kF_k^\ast\Phi,D_{2r}}\right),
\end{equation}
for every $u\in F^\ast_k\Omega^{0,q}(D_{\log k})$.
Since all the derivatives of the coefficients of the operator $\Box^{(q)}_{b,(k)}$ are uniformly bounded in $k$, it is
straightforward to see that $C_{r,\,s,k}$ can be taken to be independent of $k$ and the point $p$.
\end{proof}

\begin{theorem}\label{ss}
There is a constant $C_0>0$ such that for all $k$ and all $x\in X$, we have
\begin{equation}
k^{-n}\Pi^{q}_{\leq k\delta}(x)\leq C_0.
\end{equation}
\end{theorem}

\begin{proof}
For any $x\in X$, we choose canonical local coordinates $(z, \theta)$ defined in Theorem~ \ref{j} in a neighborhood $D$ centered at $x$. Let $s$ be a local rigid CR frame of $L$ over $D$.
For $u_k\in\mathcal H^q_{b, \leq k\delta}(X,L^k), \|u_k\|_{h^{L^k}}=1, u_k=\tilde u_k\otimes s^k$ on $D.$ Set $\tilde u_{(k)}:=k^{-\frac{n}{2}}F_k^\ast(\tilde u_k)$ on $D_{\log k}$. Write
\begin{equation}\label{r}
u_{k}=\sum\limits_{|m|\leq k\delta, m\in\mathbb Z}u_{k,m}, Tu_{k,m}=imu_{k,m}
\end{equation}
which implies that
\begin{equation}\label{w}
\tilde u_k=\sum\limits_{|m|\leq k\delta, m\in\mathbb Z}\tilde u_{k,m}, T\tilde u_{k,m}=im\tilde u_{k,m}.
\end{equation}
Then the scaling of $\tilde u_k$  given by
\begin{equation}
\tilde u_{(k)}=k^{-\frac{n}{2}}\sum\limits_{|m|\leq k\delta, m\in\mathbb Z}F_k^\ast(\Td u_{k,m})
\end{equation}
satisfies
\begin{equation}\label{u}
\|\tilde u_{(k)}\|^2_{kF_k^\ast\Phi, D_{\log k}}=\|\tilde u_k\|^2_{k\Phi, F_k(D_{\log k})}\leq\|u_k\|^2_{h^{L^k}}=1.
\end{equation}
From (\ref{s}) we have
\begin{equation}\label{v}
\Box^{(q)}_{b,(k)}\tilde u_{(k)}=0~\text{on}~D_{\log k}.
\end{equation}
By Proposition \ref{t} and combining (\ref{u}), (\ref{v}) we have
\begin{equation}\label{x}
\|\tilde u_{(k)}\|^2_{kF_k^\ast\Phi, s+1, D_r}\leq C_{r,s}(1+\|(\frac{\partial}{\partial\theta})^{s+1}\tilde u_{(k)}\|_{kF_k^\ast\Phi, D_{2r}}^2)
\end{equation}
for any $r>0$ with $D_{2r}\Subset D_{\log k}$.
Since
\begin{equation}
\begin{split}
\frac{\partial}{\partial\theta}\tilde u_{(k)}&=k^{-\frac{n}{2}}\frac{\partial}{\partial\theta}
\sum\limits_{|m|\leq k\delta, m\in\mathbb Z} F_k^\ast(\tilde u_{k,m})
=k^{-\frac{n}{2}}\sum\limits_{|m|\leq k\delta, m\in\mathbb Z}\left(\frac{im}{k}\right)\tilde u_{k,m}\left(\frac{z}{\sqrt k}, \frac{\theta}{k}\right),
\end{split}
\end{equation}
then
\begin{equation}
\left(\frac{\partial}{\partial\theta}\right)^{s+1}\tilde u_{(k)}=k^{-\frac{n}{2}}\sum\limits_{|m|\leq k\delta, m\in\mathbb Z}\left(\frac{im}{k}\right)^{s+1}\tilde u_{k,m}\left(\frac{z}{\sqrt k}, \frac{\theta}{k}\right).
\end{equation}
Thus,
\begin{equation}\label{e-ftaI}
\begin{split}
\left\|\left(\frac{\partial}{\partial\theta}\right)^{s+1}\tilde u_{(k)}\right\|^2_{kF_k^\ast\Phi, D_r}
\leq \delta^{s+1}(k\delta)\sum\limits_{|m|\leq k\delta, m\in\mathbb Z}\left\|k^{-\frac{n}{2}}\tilde u_{k,m}\left(\frac{z}{\sqrt k}, \frac{\theta}{k}\right)\right\|^2_{kF_k^\ast\Phi, D_{2r}}.
\end{split}
\end{equation}
From (\ref{w}), there is a function $\hat{\tilde u}_{k,m}(z)\in C^\infty$ such that
\begin{equation}
\tilde u_{k,m}(z,\theta)=\hat{\tilde u}_{k,m}(z)e^{im\theta}~\text{on}~D.
\end{equation}
Since
\begin{equation}\label{e-fta}
\begin{split}
&k\sum\limits_{|m|\leq k\delta, m\in\mathbb Z}\left\|k^{-\frac{n}{2}}\tilde u_{k,m}\left(\frac{z}{\sqrt k}, \frac{\theta}{k}\right)\right\|^2_{kF_k^\ast\Phi, D_{2r}}\\
&\leq\sum\limits_{|m|\leq k\delta}\int_{D_{2r}}k^{-(n-1)}\left|\hat{\tilde u}_{k,m}\left(\frac{z}{\sqrt k}\right)\right|^2e^{-k\Phi(\frac{z}{\sqrt k})}m\left(\frac{z}{\sqrt k}\right)dv(z)d\theta\\
&\leq\sum\limits_{|m|\leq k\delta}(4r)\int_{|z|\leq \frac{2r}{\sqrt k}}|\hat{\tilde u}_{k,m}(z)|^2e^{-k\Phi(z)}m(z)dv(z)\\
&\leq\frac{4r}{\varepsilon}\sum\limits_{|m|\leq k\delta}\int_{|z|\leq \frac{2r}{\sqrt k}}\int_{ |\theta|<\varepsilon}|\tilde u_{k,m}(z,\theta)|^2e^{-k\Phi(z)}m(z)dv(z)d\theta\\
&\leq \frac{4r}{\varepsilon}\sum\limits_{|m|\leq k\delta}\| u_{k,m}\|^2_{h^{L^k}}
\leq \frac{4r}{\varepsilon}\|u_{k}\|^2_{h^{L^k}}\leq\frac{4r}{\varepsilon},
\end{split}
\end{equation}
where $\varepsilon>0$ is a small constant. From (\ref{e-fta}) and (\ref{e-ftaI}), we deduce that
\begin{equation*}
\left\|\left(\frac{\partial}{\partial\theta}\right)^{s+1}\tilde u_{(k)}\right\|^2_{kF_k^\ast\Phi, D_r}\leq \tilde C_{r,s},
\end{equation*}
where $\tilde C_{r, s}$ is a constant independent of $k$. Combining this with (\ref{x}), there exists a constant $C_{r,s}^\prime>0$ independent of $k$ such that
\begin{equation}\label{y}
\|\tilde u_{(k)}\|^2_{kF_k^\ast\Phi, s+1, D_r}\leq C_{r,s}^\prime.
\end{equation}
From (\ref{y}) and Sobolev embedding theorem, there exists a constant $C(x)>0$ such that for all $k$, we have $k^{-n}|u_k(x)|_{h^{L^k}}^2=|\tilde u_{(k)}(0)|^2\leq C(x).$
Since $X$ is compact, we infer that $C^\prime=\sup\{k^{-n}|u_k(x)|^2_{h^{L^k}}:k>0, x\in X\}<\infty.$ Thus, for a local orthonormal frame $\{e^J:|J|=q, J~\text{strictly increasing}\}$ we have that $\sup\{k^{-n}S^{q}_{\leq k\delta,J}(x): x\in X, k>0\}\leq C_0$. From Lemma \ref{z}, we get the conclusion of Theorem \ref{ss}.
\end{proof}

\subsection{The Heisenberg guoup $ H_n$}

We identify $\mathbb R^{2n-1}$ with the Heisenberg gruop $H_n:=\mathbb C^{n-1}\times\mathbb R.$ We also write $(z,\theta)$ to denote the coordinates of $H_n$, $z=(z_1,\cdots, z_{n-1}),\theta\in\mathbb R,$ $z_j=x_{2j-1}+i x_{2j}, j=1,\cdots, n-1$. Then
\begin{equation}
\left\{U_{j,  H_n}=\frac{\partial}{\partial z_j}+i\lambda_j\overline z_j\frac{\partial}{\partial\theta};j=1,\cdots, n-1\right\}
\end{equation}
and
$$\left\{U_{j,  H_n}, \overline {U_{j,  H_n}}, T=\frac{\partial}{\partial\theta};j=1,\cdots,n-1\right\}$$
are local frames for the bundles of $T^{1,0}H_n$ and $\mathbb CTH_n$. Then
\begin{equation}
\left\{dz_j, d\overline z_j, \omega_0=-d\theta+\sum_{j=1}^{n-1}(i\lambda_j\overline z_jdz_j-i\lambda_jz_jd\overline z_j):j=1,\cdots, n-1\right\}
\end{equation}
is the basis of $\mathbb CT^\ast H_n$ which are dual to $\{U_{j, H_n}, \overline{U_{j, H_n}}, -T\}$. Let $\langle\,\cdot\,|\,\cdot\,\rangle$ be the Hermitian metric defined on $T^{\ast0,q}H_n$ such that $\{d\overline z^J: |J|=q; J~\text{strictly~increasing}\}$ is an orthonormal frame of $T^{\ast0, q} H_n$. Let
\begin{equation}
\overline\partial_{b, H_n}=\sum_{j=1}^{n-1}d\overline z_j\wedge\overline{U_{j, H_n}}:\Omega^{0,q}(H_n)\rightarrow\Omega^{0,q+1}(H_n)
\end{equation}
be the Cauchy-Riemann operator defined on $H_n$.
Put $\Phi_0(z)=\sum\limits_{j,t=1}^{n-1}\mu_{j,t}\overline z_jz_t\in C^\infty(H_n, \mathbb R)$. Let $(\cdot|\cdot)_{\Phi_0}$ be the inner product on $\Omega_0^{0,q}(H_n)$ with respect to the weight function $\Phi_0(z)$ defined as follows:
\begin{equation}
(f|g)_{\Phi_0}=\int_{H_n}\langle f|g\rangle e^{-\Phi_0(z)}dv(z)d\theta
\end{equation}
where $dv(z)=2^{n-1}dx_1\cdots dx_{2n-2}$. We denote by $\|\cdot\|_{\Phi_0}$ the norm on $\Omega^{0, q}_0(H_n, \Phi_0)$ induced by the inner product $(\cdot|\cdot)_{\Phi_0}$. Let us denote by $L^2_{(0, q)}(H_n, \Phi_0)$ the completion of $\Omega^{0, q}_0(H_n)$ with respect to the norm $\|\cdot\|_{\Phi_0}$. Let $\overline\partial_{b, H_n}^{\ast,\Phi_0}: \Omega^{0,q+1}(H_n)\rightarrow\Omega^{0,q}(H_n)$ be the formal adjoint of $\overline\partial_{b, H_n}$ respect to $(\cdot|\cdot)_{\Phi_0}$. The Kohn Laplacian on $H_n$ is given by
\begin{equation}
\Box^{(q)}_{b, H_n}=\overline\partial_{b, H_n}\overline\partial_{b, H_n}^{\ast,\Phi_0}+\overline\partial_{b, H_n}^{\ast,\Phi_0}\overline\partial_{b, H_n}:\Omega^{0,q}(H_n)\rightarrow\Omega^{0,q}(H_n).
\end{equation}

We pause and introduce some notations.
Choose $\chi(\theta)\in C_0^\infty(\mathbb R)$
so that $\chi(\theta)=1$ when $|\theta|<1$ and
$\chi(\theta)=0$ when $|\theta|>2$
and set $\chi_j(\theta)=\chi(\frac{\theta}{j}),j\in\mathbb N$. For any $u(z, \theta)\in\Omega^{0, q}(H_n)$ with $\|u\|_{\Phi_0}<\infty$. Let
\begin{equation}
\hat u_j(z, \eta)=\int_{\mathbb R} u(z, \theta)\chi_j(\theta)e^{-i\theta\eta}d\theta\in\Omega^{0, q}(H_n), j=1, 2,\ldots.
\end{equation}
From Parseval's formula, $\{\hat u_j(z, \eta)\}$ is a cauchy sequence in $L^2_{(0, q)}(H_n, \Phi_0)$. Thus there is $\hat u(z, \eta)\in L^2_{(0, q)}(H_n, \Phi_0)$ such that $\hat u_j(z, \eta)\rightarrow \hat u(z, \eta)$ in $L^2_{(0, q)}(H_n, \Phi_0)$. We call $\hat u(z, \eta)$ the partial Fourier transform of $u(z, \theta)$ with respect to $\theta$. From Parseval's formula, we can check that
\begin{equation} \label{s4-e12-1}
\begin{split}
\int_{H_n}\!\abs{\hat u(z, \eta)}^2e^{-\Phi_0(z)}dv(z)d\eta=
2\pi\int_{H_n}\!\abs{u(z, \theta)}^2e^{-\Phi_0(z)}dv(z)d\theta.
\end{split}
\end{equation}
Let $s\in L^2_{(0,q)}(H_n, \Phi_0)$.
Assume that $\int\!\abs{s(z, \eta)}^2d\eta<\infty$ and $\int\!\abs{s(z, \eta)}d\eta<\infty$ for all $z\in\Complex^{n-1}$. Then, from Parseval's formula, we can check that
\begin{equation} \label{s4-e12-2}
\begin{split}
&\iint\!\langle\,\hat u(z, \eta)\,|\,s(z, \eta)\,\rangle e^{-\Phi_0(z)}d\eta dv(z)\\
&=\iint\!\langle\,u(z, \theta)\,|\,\int\! e^{i\theta\eta}s(z, \eta)d\eta\,\rangle e^{-\Phi_0(z)}d\theta dv(z).
\end{split}
\end{equation}

\subsection{Proof of Theorem~\ref{m}}

Now we can prove the second part of Theorem \ref{m}.

\begin{theorem}\label{eee}
\begin{equation}
\limsup_{k\rightarrow\infty}k^{-n}\Pi^{q}_{\leq k\delta}(x)\leq(2\pi)^{-n}\int_{\mathbb R_{x,q}\bigcap[-\delta,\delta]}|\det(\mathcal R_x^L+2s\mathcal L_x)|ds
\end{equation}
for all $x\in X$.
\end{theorem}

\begin{proof}
Fix $x\in X$ and let $s$ be a rigid CR frame of $L$ on an open neighborhood $D$ of $x$ and $|s|^2_{h^L}=e^{-\Phi}$. We take  canonical local coordinates $(z,\theta), z_j=x_{2j-1}+i x_{2j}, j=1,\cdots, n-1$ defined in Theorem \ref{j} such that $\omega_0(x)=-d\theta, (z(x), \theta(x))=0$, and \eqref{e-gue150303}, \eqref{e-gue150303I}, \eqref{e-gue150303II} hold . Until further notice, we work with canonical coordinates $(z,\theta)$ defined on an open neighbourhood $D$ of $x$ and we identify $D$ with some open set in $\mathbb R^{2n-1}$. We will use the same notations as in section~\ref{kk}. Fix $\abs{J}=q$, $J$ is strictly increasing.
First, from definition of extremal function, there exists a sequence $\alpha_{k_j}\in\mathcal H^q_{b, \leq k_j\delta}(X, L^{k_j})$, $0<k_1<k_2<\cdots$, such that
$\|\alpha_{k_j}\|^2_{h^{L^{k_j}}}=1$ and
\begin{equation}\label{bbb}
\lim_{j\rightarrow\infty}k_j^{-n}|\alpha_{k_j, J}(x)|^2_{h^{L^{k_j}}}=\limsup_{k\rightarrow\infty}k^{-n} S^{q}_{\leq k\delta, J}(x)
\end{equation}
where $\alpha_{k_j, J}$ is the component of $\alpha_{k_j}$ along the direction $e^J$. Put
$\alpha_{k_j}=\tilde{\alpha}_{k_j}\otimes s^{k_j}$, $\tilde{\alpha}_{k_j}\in\Omega^{0, q}(D)$. We will always use  $\alpha_{k_j}$ to denote $\tilde\alpha_{k_j}$ if there is no misunderstanding, then
\begin{equation}
\alpha_{(k_j)}=k_j^{-\frac{n}{2}}F_{k_j}^\ast(\alpha_{k_j})\in F_{k_j}^\ast\Omega^{0,q}(D_{\log k_j}).
\end{equation}
It is easy to see that for every $j$,
\begin{equation}\label{ddd}
\|\alpha_{(k_j)}\|_{k_jF_{k_j}^\ast\Phi, D_{\log k_j}}\leq 1,\ \ \Box^{(q)}_{b,(k_j)}\alpha_{(k_j)}=0.
\end{equation}
Moreover, we can repeat the procedure in the proof of Theorem~\ref{ss} and obtain that for every $r>0$ and $s\in\mathbb N_0$ there is a $C_{r,s}$ independent of $k_j$ such that
\begin{equation}\label{e-fb}
\norm{(\frac{\pr}{\pr\theta})^{s+1}\alpha_{(k_j)}}_{k_jF^*_{k_j}\Phi,D_r}\leq C_{r,s},\ \ \forall j.
\end{equation}
From \eqref{ddd}, \eqref{e-fb} and Proposition~\ref{t}, we can repeat the same argument in Theorem 2.9 of \cite{HM12} and conclude that there is a subsequence $\{\alpha_{(k_{s_1})},\alpha_{(k_{s_2})},\ldots\}$ of $\{\alpha_{(k_j)}\}$, $0<k_{s_1}<k_{s_2}<\cdots$, such that $\alpha_{(k_{s_t})}$ converges uniformly with all derivatives on any compact subset of $H_n$ to a smooth form $u=\sum\nolimits_{|J|=q}^\prime u_Jd\overline z^J\in\Omega^{0,q}(H_n)$ as $t\To\infty$. Thus,
\begin{equation}\label{e-gue150303III}
\limsup_{k\To\infty}k^{-n}S^{q}_{\leq k\delta,J}(x)\leq\abs{u_J(0)}^2.
\end{equation}
Moreover, (\ref{ddd}) implies that $u$ satisfies
\begin{equation}
\|u\|_{\Phi_0}\leq 1,\ \  \Box^{(q)}_{b, H_n}u=0.
\end{equation}
Then we will need

\begin{lemma}\label{ee}
With the notations above, $\hat u(z,\eta)\equiv0$ in $L^2_{(0,q)}(H_n,\Phi_0)$ when  $|\eta|>\delta$.
\end{lemma}

\begin{proof}
To prove $\hat u(z,\eta)\equiv0$ when $|\eta|>\delta$, we only need to show that for any $\varphi(z,\eta)\in C_0^\infty(\mathbb C^{n-1}\times\{\eta\in\mathbb R: |\eta|>\delta\})$ and $\abs{J}=q$, $J$ is strictly increasing, we have
\begin{equation}
\int_{H_n}\hat u_J(z,\eta)\varphi(z,\eta)e^{-\Phi_0(z)}dv(z)d\eta=0.
\end{equation}
We assume that ${\rm supp\,}\varphi\Subset \{z\in\mathbb C^{n-1}:|z|\leq r_0\}\times\{\eta\in\mathbb R: |\eta|>\delta\}$. Here, $\{z\in\mathbb C^{n-1}:|z|<r_0\}$ means that $\{z\in\mathbb C^{n-1}: |z_j|<r_0, j=1,\cdots, n-1\}.$
Choose $\chi\in C_0^\infty(\mathbb R)$ such that $\chi\equiv1$ when $|\theta|\leq 1$ and ${\rm supp\,}\chi\Subset\{\theta\in\mathbb R: |\theta|<2\}$. From \eqref{s4-e12-2}, we have
\begin{equation}
\begin{split}
\frac{1}{2\pi}\int_{H_n}\hat u_J(z,\eta)\varphi(z,\eta)e^{-\Phi_0(z)}dv(z)d\eta
&=\int_{H_n}u_J(z,\theta)\check{\varphi}(z,\theta)e^{-\Phi_0(z)}dv(z)d\theta\\
&=\lim_{r\rightarrow\infty}\int_{H_n}u_J(z,\theta)\check{\varphi}(z,\theta)e^{-\Phi_0(z)}
\chi(\frac{\theta}{r})dv(z)d\theta,
\end{split}
\end{equation}
where $\check{\varphi}(z,\theta):=\frac{1}{2\pi}\int_{\mathbb R} e^{i\theta\eta}\varphi(z,\eta)d\eta$ is the inverse Fourier transform of $\varphi(z,\eta)$ respect to $\eta$. For simplicity, we may assume that $\alpha_{(k_{j})}$ converges uniformly with all derivatives on any compact subset of $H_n$ to $u$ as $j\To\infty$. Note that $\alpha_{k_j}\in\Omega^{0,q}_{\leq k_j\delta}(X,L^k)$. For each $j$, on $D$, we can write
\[\alpha_{k_j}=s^{k_j}\otimes\sum_{m\in\mathbb Z,\abs{m}\leq k_j\delta}\Td\alpha_{k_j,m},\ \ \Td\alpha_{k_j,m}\in\Omega^{0,q}(D),\ \ \forall m\in\mathbb Z,\ \abs{m}\leq k_j\delta,\]
and for each $m\in\mathbb Z$, $\abs{m}\leq k_j\delta$, we can write
\[\begin{split}
&\Td\alpha_{k_j,m}=\sideset{}{'}\sum_{\abs{J}=q}\Td\alpha_{k_j,m,J}(z,\theta)e^J(z),\\
&\mbox{$\Td\alpha_{k_j,m,J}=\hat\alpha_{k_j,m,J}(z)e^{im\theta}$, $\hat\alpha_{k_j,m,J}(z)\in C^\infty(D)$, $\forall\abs{J}=q$, $J$ is strictly increasing}.
\end{split}\]
When $r$ is fixed, by dominated convergence theorem
\begin{equation}\label{e-gue150303b}
\begin{split}
~~~~&\int_{H_n}u_J(z,\theta)\check{\varphi}(z,\theta)e^{-\Phi_0(z)}
\chi(\frac{\theta}{r})dv(z)d\theta\\
&=\lim_{j\rightarrow\infty}\sum_{|m|\leq k_j\delta}\int_{H_n}k_j^{-\frac{n}{2}}
\hat\alpha_{k_j,m,J}\left(\frac{z}{\sqrt {k_j}}\right)e^{i\frac{m}{ k_j}\theta}\check\varphi(z,\theta)\chi\left(\frac{\theta}{r}\right)e^{-\Phi_0(z)}dv(z)d\theta\\
&=\lim_{j\rightarrow\infty}\sum_{|m|\leq k_j\delta}\int_{|z|\leq r_0}\int_{\mathbb R}
k_j^{-\frac{n}{2}}\hat\alpha_{k_j,m,J}\left(\frac{z}{\sqrt {k_j}}\right)e^{i\frac{m}{k_j}\theta}
\check\varphi(z,\theta)\chi\left(\frac{\theta}{r}\right)e^{-\Phi_0(z)}dv(z)d\theta.
\end{split}
\end{equation}
Since ${\rm supp\,}\varphi(z,\eta)\Subset\{z\in\mathbb C^{n-1}:|z|\leq r_0\}\times\{\theta\in\mathbb R: |\eta|>\delta\}$ and $|\frac{m}{k_j}|<\delta$, we have
\begin{equation}\label{ff}
\begin{split}
&\sum_{|m|\leq k_j\delta}\int_{H_n}k_j^{-\frac{n}{2}}\hat\alpha_{k_j, m, J}\left(\frac{z}{\sqrt{k_j}}\right)e^{i\frac{m}{k_j}\theta}\check\varphi(z,\theta)e^{-\Phi_0(z)}dv(z)d\theta\\
&=\sum_{|m|\leq k_j\delta}\int_{H_n}k_j^{-\frac{n}{2}}\hat\alpha_{k_j, m, J}\left(\frac{z}{\sqrt{k_j}}\right)\varphi(z,-\frac{m}{k_j})e^{-\Phi_0(z)}dv(z)=0.
\end{split}
\end{equation}
By (\ref{ff})
\begin{equation}\label{e-gue150303bI}
\begin{split}
&\lim_{j\rightarrow\infty}\sum_{|m|\leq k_j\delta}\int_{|z|\leq r_0}\int_{\mathbb R}
k_j^{-\frac{n}{2}}\hat\alpha_{k_j,m,J}\left(\frac{z}{\sqrt {k_j}}\right)e^{i\frac{m}{k_j}\theta}
\check\varphi(z,\theta)\chi\left(\frac{\theta}{r}\right)e^{-\Phi_0(z)}dv(z)d\theta\\
&=\lim_{j\rightarrow\infty}\sum_{|m|\leq k_j\delta}\int_{|z|\leq r_0}\int_{\mathbb R}k_j^{-\frac{n}{2}}\hat\alpha_{k_j,m, J}\left(\frac{z}{\sqrt{k_j}}\right)e^{i\frac{m}{k_j}\theta}
\check\varphi(z,\theta)\left(\chi\left(\frac{\theta}{r}\right)-1\right)e^{-\Phi_0(z)}dv(z)d\theta.
\end{split}
\end{equation}
Now,
\begin{equation}\label{e-gue150303bII}
\begin{split}
&\left|\sum_{|m|\leq k_j\delta}\int_{|z|\leq r_0}\int_{\mathbb
R}k_j^{-\frac{n}{2}}\hat\alpha_{k_j,m, J}\left(\frac{z}{\sqrt{k_j}}\right)e^{i\frac{m}{k_j}\theta}
\check\varphi(z,\theta)\left(\chi\left(\frac{\theta}{r}\right)-1\right)e^{-\Phi_0(z)}dv(z)d\theta\right|\\
&\leq\sum_{|m|\leq k_j\delta}\int_{|z|\leq r_0}\int_{|\theta|\geq
r}k_j^{-\frac{n}{2}}\left|\hat\alpha_{k_j,m, J}\left(\frac{z}{\sqrt{k_j}}\right)\right|\cdot|\check\varphi(z,\theta)|e^{-\Phi_0(z)}dv(z)d\theta.
\end{split}
\end{equation}
By H\"older inequality,  we have
\begin{equation}\label{bb}
\begin{split}
&\int_{|z|\leq r_0}\int_{|\theta|\geq
r}k_j^{-\frac{n}{2}}|\hat\alpha_{k_j,m, J}\left(\frac{z}{\sqrt{k_j}}\right)|\cdot|\check\varphi(z,\theta)|e^{-\Phi_0(z)}dv(z)d\theta\\
&\leq \left(\int_{|z|\leq r_0}\int_{|\theta|\geq
r}k_j^{-n}|\hat\alpha_{k_j, m,
J}\left(\frac{z}{\sqrt{k_j}}\right)|^2\cdot|\check\varphi(z,\theta)|e^{-\Phi_0(z)}dv(z)d\theta
\right)^{\frac12}\times\\
&\left(\int_{|z|\leq r_0}\int_{|\theta|\geq r}|\check\varphi(z,\theta)|e^{-\Phi_0(z)}dv(z)d\theta\right)^{\frac12}.
\end{split}
\end{equation}
Since ${\rm supp\,}\varphi(z,\eta)\Subset{\{z\in\mathbb C^{n-1}:|z|\leq r_0}\}\times\mathbb R$, we have
\begin{equation}\label{cc}
\sup_{|z|\leq r_0}|\check\varphi(z,\theta)|\leq C_{r_0}\frac{1}{|\theta|^p}, \forall ~|\theta|>>1
\end{equation}
for some $p>3, p\in\mathbb Z$ and constant $C_{r_0}>0.$
Combining (\ref{bb}) and (\ref{cc}),we have
\begin{equation}\label{dd}
\begin{split}
&\int_{|z|\leq r_0}\int_{|\theta|\geq
r}k_j^{-\frac{n}{2}}\left|\hat\alpha_{k_j,m, J}\left(\frac{z}{\sqrt{k_j}}\right)\right|\cdot|\check\varphi(z,\theta)|e^{-\Phi_0(z)}dv(z)d\theta\\
&\leq C_{r_0}\frac{1}{r^{p-1}}\left(\int_{|z|\leq r_0}k_j^{-n}\left|\hat\alpha_{k_j, m, J}\left(\frac{z}{\sqrt{k_j}}\right)\right|^2e^{-\Phi_0(z)}dv(z)\right)^{\frac12}\\
&\leq C_{r_0}^\prime\frac{1}{r^{p-1}}\left(\int_{|z|\leq r_0}k_j^{-n}\left|\hat\alpha_{k_j,m, J}\left(\frac{z}{\sqrt{k_j}}\right)\right|^2e^{-k_jF_{k_j}^\ast\Phi(z)}(F_{k_j}^\ast m)dv(z)\right)^{\frac12}\\
&\leq C_{r_0}^\prime\frac{1}{r^{p-1}}\left(\int_{|z|\leq\frac{r_0}{\sqrt{k_j}}}\frac{1}{k_j}
|\hat\alpha_{k_j,m,J}(z)|^2e^{-k_j\Phi(z)}m(z)dv(z)\right)^{\frac12},
\end{split}
\end{equation}
 for $r>>1$, where $C_{r_0}>0$ and $C_{r_0}^\prime>0$ are constants.
Then from (\ref{dd}) and Cauchy-Schwartz inequality,
\begin{equation}\label{e-gue150303bIII}
\begin{split}
~~~~~~~~&\sum_{|m|\leq k_j\delta}\int_{|z|\leq r_0}\int_{|\theta|\geq
r}k_j^{-\frac{n}{2}}\left|\hat\alpha_{k_j,m, J}\left(\frac{z}{\sqrt{k_j}}\right)\right|\cdot|\check\varphi(z,\theta)|e^{-\Phi_0(z)}dv(z)d\theta\\
&\leq C_{r_0}^\prime\frac{1}{r^{p-1}} \sqrt{\delta}\left( \sum_{|m|\leq k_j\delta}\int_{|z|\leq\frac{r_0}{\sqrt{k_j}}}|\hat\alpha_{k_j,m,J}(z)|^2
e^{-k_j\Phi(z)}m(z)dv(z)\right)^{\frac12}\\
&\leq C_{r_0}^\prime\frac{1}{r^{p-1}} \frac{\sqrt{\delta}}{\sqrt{2\varepsilon}}
\left(\sum_{|m|\leq k_j\delta}\int_{|z|\leq\frac{r_0}{\sqrt{k_j}}, |\theta|\leq\varepsilon}|\Td\alpha_{k_j,m,J}(z,\theta)|^2
e^{-k_j\Phi(z)}m(z)dv(z)d\theta\right)^{\frac12}\\
&\leq C_{r_0}^\prime\frac{1}{r^{p-1}} \frac{\sqrt{\delta}}{\sqrt{2\varepsilon}}\left(\sum_{|m|\leq k_j\delta}\|\alpha_{k_j, m}\|^2_{h^{L^{k_j}}}\right)^{\frac{1}{2}}
\leq C_{r_0}^\prime\frac{1}{r^{p-1}} \frac{\sqrt{\delta}}{\sqrt{2\varepsilon}}\|\alpha_{k_j}\|^2_{h^{L^{k_j}}}
\leq C_{r_0}^\prime\frac{1}{r^{p-1}} \frac{\sqrt{\delta}}{\sqrt{2\varepsilon}}.
\end{split}
\end{equation}
From \eqref{e-gue150303b}, \eqref{e-gue150303bI}, \eqref{e-gue150303bII} and \eqref{e-gue150303bIII}, we get
\begin{equation}
\begin{split}
\left|\int_{H_n}u_J(z,\theta)\check{\varphi}(z,\theta)e^{-\Phi_0(z)}
\chi(\frac{\theta}{r})dv(z)d\theta\right|
\leq C_{r_0}^\prime\frac{1}{r^{p-1}} \frac{\sqrt{\delta}}{\sqrt{2\varepsilon}}.
\end{split}
\end{equation}
Letting $r\rightarrow\infty$, we get the conclusion of Lemma \ref{ee}.
\end{proof}

We pause and introduce some notations.
For fixed $\eta\in\mathbb R$, put $\Phi_{\eta}(z)=-2\eta\sum\limits_{j=1}^{n-1}\lambda_j|z_j|^2+\sum\limits_{j,t=1}^{n-1}\mu_{j,t}\overline z_jz_t.$ We take the Hermitian metric $\langle\,\cdot\,|\,\cdot\,\rangle$ on the bundle $T^{\ast0,q}\Complex^{n-1}$ of $(0, q)$ forms on $\Complex^{n-1}$ so that $\{d\ol z^J; \text{$\abs{J}=q$, $J$ strictly increasing}\}$ is an orthonormal frame. We also let $\Omega^{0,q}(\Complex^{n-1})$ denote the space of smooth sections of $T^{\ast0,q}\Complex^{n-1}$ over $\Complex^{n-1}$ and let $\Omega^{0,q}_0(\Complex^{n-1})$ be the subspace of $\Omega^{0,q}(\Complex^{n-1})$ whose elements have compact support in $\Complex^{n-1}$. Let $(\,\cdot\,|\,\cdot\, )_{\Phi_\eta}$ be the inner product on $\Omega^{0,q}_0(\Complex^{n-1})$ defined by
\[(\,f\,|\,g\,)_{\Phi_\eta}=\int_{\Complex^{n-1}}\!\langle\,f\,|\,g\,\rangle e^{-\Phi_\eta(z)}dv(z)\,,\quad f, g\in\Omega^{0,q}_0(\Complex^{n-1})\,,\] where $dv(z)=2^{n-1}dx_1dx_2\cdots dx_{2n-2}$, and let $\norm{\cdot}_{\Phi_\eta}$ denote the corresponding norm. Let us denote by $L^2_{(0,q)}(\mathbb C^{n-1}, \Phi_\eta)$ the completion of $\Omega_0^{0,q}(\mathbb C^{n-1})$ with respect to the norm $\|\cdot\|_{\Phi_\eta}$. Let
$$\Box^{(q)}_{\Phi_\eta}=\overline\partial^{\ast,\Phi_\eta}\overline\partial+\overline\partial\,
\overline\partial^{\ast,\Phi_\eta}:\Omega^{0,q}(\mathbb C^{n-1})\rightarrow\Omega^{0,q}(\mathbb C^{n-1})$$
be the complex Laplacian with respect to $(\,\cdot\,|\,\cdot\,)_{\Phi_\eta}$, where $\ddbar^{\ast,\Phi_\eta}$ is the formal adjoint of $\ddbar$ with respect to $(\,\cdot\,|\,\cdot\,)_{\Phi_\eta}$.
Let $B^{(q)}_{\Phi_\eta}:L^2_{(0,q)}(\mathbb C^{n-1}, \Phi_\eta)\rightarrow\rm{Ker}\Box^{(q)}_{\Phi_\eta}$ be the Bergman projection and $B^{(q)}_{\Phi_{\eta}}(z,w)$ be the distribution kernel of $B^{(q)}_{\Phi_\eta}$ with respect to $(\,\cdot\,|\,\cdot\,)_{\Phi_\eta}$(see section 3.2 in \cite{HM12}). Let $M_{\Phi_\eta}: T_z^{1,0}\mathbb C^{n-1}\rightarrow T_z^{1,0}\mathbb C^{n-1}, z\in\mathbb C^{n-1}$ be the linear map defined by $$\langle M_{\Phi_\eta}U|V\rangle=\partial\overline\partial\Phi_{\eta}(U, \overline V), U, V\in T_z^{1,0}\mathbb C^{n-1}.$$
Put $$\mathbb R_q=\{\eta\in\mathbb R: M_{\Phi_{\eta}}~\text{ has exactly $q$ negative eigenvalues and $n-1-q$ positive eigenvalues}\}.$$
The following lemma is well known(see Berman \cite{Be04}, Hsiao and Marinescu \cite{HM12}, Ma and Marinescu \cite{MM07} ).

\begin{lemma}\label{abd}
If $\eta\not\in\mathbb R_q$, then $B^{(q)}_{\Phi_\eta}(z, z)=0$ for all $z\in\mathbb C^{n-1}.$ If $\eta\in\mathbb R_q$, then
\begin{equation}
\sum\nolimits_{|J|=q}^\prime\langle B^{(q)}_{\Phi_\eta}(z,z)d\overline z^J|d\overline z^J\rangle=e^{\Phi_\eta(z)}(2\pi)^{-n+1}|\det M_{\Phi_\eta}|\cdot 1_{\mathbb R_q}(\eta)
\end{equation}
\end{lemma}
The following is also well-known.
\begin{lemma}[See Theorem 3.1 and Lemma 3.5 in \cite{HM12}]\label{ii}
For almost all $\eta\in\Real$, $\hat u(z,\eta)\in\Omega^{0,q}(\Complex^{n-1})$, $\int_{\Complex^{n-1}}\abs{\hat u(z,\eta)}^2e^{-\Phi_0(z)}dv(z)<\infty$ and
\[
|\hat u_J(z,\eta)|^2\leq \exp{(\eta\sum_{j=1}^{n-1}\lambda_j|z_j|^2)}\langle B^{(q)}_{\Phi_{\eta}}(z,z)d\overline z^J|d\overline z^J\rangle\int_{\mathbb C^{n-1}}|\hat u(w,\eta)|^2e^{-\Phi_0(w)}dv(w),
\]
for every strictly increasing index $J$, $\abs{J}=q$.
\end{lemma}

Now, we can prove

\begin{proposition}\label{p-gue150306}
For every $|J|=q$, $J$ is strictly increasing, we have
\[u_J(0,0)=\frac{1}{2\pi}\int_{\eta\in\mathbb R,\abs{\eta}\leq\delta}\hat u_J(0,\eta)d\eta.\]
\end{proposition}

\begin{proof}
Let $\chi\in C^\infty_0(\Real)$, $\int_\Real\!\chi d\theta=1$, $\chi\geqslant0$ and $\chi_\varepsilon\in C^\infty_0(\Real)$, $\chi_\varepsilon(\theta)=\frac{1}{\varepsilon}\chi(\frac{\theta}{\varepsilon})$. Then,
$\chi_\varepsilon\To\delta_0$, $\varepsilon\To0^+$
in the sense of distributions.
Let $\hat\chi_\varepsilon:=\int e^{-i\theta\eta}\chi_\varepsilon(\theta)d\theta$
be the Fourier transform of $\chi_\varepsilon$. We can check that
$\abs{\hat\chi_\varepsilon(\eta)}\leqslant  1$ for all $\eta\in\Real$, $\hat\chi_\varepsilon(\eta)=\hat\chi(\varepsilon\eta)$ and
$\lim_{\varepsilon\To0}\hat\chi_\varepsilon(\eta)=\lim_{\varepsilon\To0}\hat\chi(\varepsilon\eta)=\hat\chi(0)=1$. Let $\varphi\in C^\infty_0(\Complex^{n-1})$ such that $\int_{\Complex^{n-1}}\!\varphi(z)dv(z)=1$, $\varphi\geqslant0$, $\varphi(z)=0$ if $\abs{z}>1$. Put
$g_j(z)=j^{2n-2}\varphi(jz)e^{\Phi_0(z)}$, $j=1,2,\ldots$. Then, for $J$ is strictly increasing, $\abs{J}=q$, we have
\begin{equation} \label{s4-e41}
u_J(0,0)=\lim_{j\To\infty}\lim_{\varepsilon\To0^+}\int_{H_n}\!\langle\,u(z,\theta)\,|\,\chi_\varepsilon(\theta)g_j(z)d\ol z^J\,\rangle e^{-\Phi_0(z)}dv(z)d\theta.
\end{equation}
From \eqref{s4-e12-2}, we see that
\begin{equation} \label{s4-e42}
\begin{split}
&\iint\!\langle\,u(z,\theta)\,|\,\chi_\varepsilon(\theta)g_j(z)d\ol z^J\,\rangle e^{-\Phi_0(z)}dv(z)d\theta\\
&=\frac{1}{2\pi}\iint\!\langle\hat u(z, \eta)\,|\,\hat\chi_\varepsilon(\eta)g_j(z)d\ol z^J\,\rangle e^{-\Phi_0(z)}d\eta dv(z).
\end{split}
\end{equation}
From Lemma~\ref{ee}, it is not difficult to check that for every $j$ and every $\varepsilon>0$,
\begin{equation} \label{e-gue150306}
\begin{split}
&\iint\!\langle\hat u(z, \eta)\,|\,\hat\chi_\varepsilon(\eta)g_j(z)d\ol z^J\,\rangle e^{-\Phi_0(z)}d\eta dv(z)\\
&=\iint_{\abs{\eta}\leq\delta}\langle\hat u(z, \eta)\,|\,\hat\chi_\varepsilon(\eta)g_j(z)d\ol z^J\,\rangle e^{-\Phi_0(z)}d\eta dv(z).
\end{split}
\end{equation}
From Lemma~\ref{ii}, \eqref{e-gue150306}, we can apply Lebesque dominated convergence theorem and conclude that
\begin{equation} \label{s4-e42-1}
\begin{split}
\lim_{\varepsilon\To0^+}&\iint\!\langle\,\hat u(z, \eta)\,|\,\hat\chi_\varepsilon(\eta)g_j(z)d\ol z^J\,\rangle e^{-\Phi_0(z)}d\eta dv(z)\\
&=\iint_{\abs{\eta}\leq\delta}\langle\,\hat u(z, \eta)\,|\,g_j(z)d\ol z^J\,\rangle e^{-\Phi_0(z)}d\eta dv(z).
\end{split}
\end{equation}
From \eqref{s4-e42} and \eqref{s4-e42-1}, \eqref{s4-e41} becomes
\begin{equation} \label{s4-e42-2}
u_J(0, 0)=\lim_{j\To\infty}\frac{1}{2\pi}\iint_{\abs{\eta}\leq\delta}\langle\,\hat u(z, \eta)\,|\,g_j(z)d\ol z^J\,\rangle e^{-\Phi_0(z)}d\eta dv(z).
\end{equation}
Put $f_j(\eta)=\frac{1}{2\pi}\int\!\langle\hat u(z, \eta)\,|\,g_j(z)d\ol z^J\,\rangle e^{-\Phi_0(z)}dv(z)$.
Since $\hat u(z, \eta)\in\Omega^{0,q}(\Complex^{n-1})$ for almost all $\eta$, we have
$\lim_{j\To\infty}f_j(\eta)=\frac{1}{2\pi}\hat u_J(0, \eta)$
almost everywhere. Now, for almost every $\eta\in\Real$,
\begin{equation} \label{s4-e42-3}
\begin{split}
\abs{f_j(\eta)}
&=\frac{1}{2\pi}\abs{\int_{\abs{z}\leqslant \frac{1}{j}}\!
\langle\,\hat u(z, \eta)\,|\,j^{2n-2}\varphi(jz)d\ol z^J\,\rangle dv(z)}\\
&\leqslant \frac{1}{2\pi}\Bigr(\int_{\abs{z}\leqslant \frac{1}{j}}\!
\abs{\hat u(z, \eta)}^2e^{-\Phi_0(z)}j^{2n-2}dv(z)
\Bigr)^{\frac{1}{2}}
\Bigr(\int_{\abs{z}\leqslant \frac{1}{j}}\!\abs{\varphi(jz)}^2e^{\Phi_0(z)}j^{2n-2}dv(z)\Bigr)^{\frac{1}{2}}\\
&\leqslant  C_1\Bigr(\int_{\abs{z}\leqslant 1}\!
\abs{\hat u(\frac{z}{j}, \eta)}^2e^{-\Phi_0(z/j)}dv(z)\Bigr)^{\frac{1}{2}}\\
&\leqslant  C_2\Bigr(\int_{\abs{z}\leqslant 1}\!e^{2\eta\sum^{n-1}_{t=1}\lambda_t\abs{\frac{z_t}{j}}^2}\abs{{\rm Tr\,}B^{(q)}_{\Phi_\eta}(\frac{z}{j},\frac{z}{j})}dv(z)\Bigr)^{\frac{1}{2}}\\
&\times\Bigr(\int_{\Complex^{n-1}}\abs{\hat u(w, \eta)}^2e^{-\Phi_0(w)}dv(w)\Bigr)^{\frac{1}{2}}\ \ (\mbox{here we used Lemma~\ref{ii}})\\
&\leqslant  C_3\Bigr(\int_{\Complex^{n-1}}\abs{\hat u(w, \eta)}^2e^{-\Phi_0(w)}dv(w)\Bigr)^{\frac{1}{2}},
\end{split}
\end{equation}
where $C_1, C_2, C_3$ are positive constants. From this and the Lebesgue dominated convergence theorem, we conclude that
\[
u_J(0, 0)=\lim_{j\To\infty}\int_{\abs{\eta}\leq\delta}f_j(\eta)d\eta=\int_{\abs{\eta}\leq q}\lim_{j\To\infty}f_j(\eta)d\eta
=\frac{1}{2\pi}\int\hat u_J(0, \eta)d\eta\,.
\]
\end{proof}


Now we turn to our situation. Fix $\abs{J}=q$, $J$ is strictly increasing.
By Proposition~\ref{p-gue150306}, Lemma~\ref{ii} and notice that (see \eqref{s4-e12-1})
\[\int_{|\eta|\leq\delta}|\hat u(w,\eta)|^2e^{-\Phi_0(w)}dv(w)d\eta\leq \int|\hat u(w,\eta)|^2e^{-\Phi_0(w)}dv(w)d\eta\leq 2\pi,\]
we have
\begin{equation}\label{e-gue150306ab}
\begin{split}
&\abs{u_J(0,0)}=\frac{1}{2\pi}\abs{\int_{|\eta|\leq\delta}\hat u_J(0,\eta)d\eta}\\
&\leq\frac{1}{2\pi}\int_{|\eta|\leq\delta}|\hat u_J(0,\eta)|\frac{(\int_{\mathbb C^{n-1}}|\hat u(w,\eta)|^2e^{-\Phi_0(w)}dv(w))^{\frac12}}{(\int_{\mathbb C^{n-1}}|\hat u(w,\eta)|^2e^{-\Phi_0(w)}dv(w))^{\frac12}}d\eta\\
\leq&\frac{1}{2\pi}\left(\int_{|\eta|\leq\delta}\frac{|\hat u_J(0,\eta)|^2}{\int_{\mathbb C^{n-1}}|\hat u(w,\eta)|^2e^{-\Phi_0(w)}dv(w)}d\eta\right)^{\frac12}\times\\
&\left(\int_{|\eta|\leq\delta}|\hat u(w,\eta)|^2e^{-\Phi_0(w)}dv(w)d\eta\right)^{\frac12}\\
&\leq\frac{1}{\sqrt{2\pi}}\left(\int_{|\eta|\leq\delta}\frac{\langle B^{(q)}_{\Phi_{\eta}}(0,0)d\overline z^J|d\overline z^J\rangle \int_{\mathbb C^{n-1}}|\hat u(w,\eta)|^2e^{-\Phi_0(w)}dv(w)}{\int_{\mathbb C^{n-1}}|\hat u(w,\eta)|^2e^{-\Phi_0(w)}dv(w)}d\eta\right)^{\frac12}\\
&\leq\frac{1}{\sqrt{2\pi}}\left(\int_{|\eta|\leq\delta}\langle B^{(q)}_{\Phi_\eta}(0,0)d\overline z^J|d\overline z^J\rangle d\eta\right)^{\frac12}.
\end{split}
\end{equation}
Combining \eqref{e-gue150303III} and \eqref{e-gue150306ab}, we have
\begin{equation}\label{ccc}
\limsup_{k\rightarrow\infty}k^{-n} S^{q}_{\leq k\delta, J}(p)\leq \frac{1}{2\pi}\int_{|\eta|\leq\delta}\langle B^{(q)}_{\Phi_\eta}(0,0)d\overline z^J|d\overline z^J\rangle d\eta.
\end{equation}
Then Lemma \ref{z}, Lemma~\ref{abd} and (\ref{ccc}) imply that
\begin{equation}
\begin{split}
\limsup_{k\rightarrow\infty}k^{-n}\Pi^q_{\leq k\delta}(x)
&\leq \frac{1}{2\pi}\int_{|\eta|\leq\delta}\sum\nolimits_{|J|=q}^\prime\langle B^{(q)}_{\Phi_\eta}(0,0)d\overline z^J|d\overline z^J\rangle d\eta\\
&\leq \frac{1}{(2\pi)^n}\int_{|\eta|\leq\delta}|\det M_{\Phi_{\eta}}|\cdot 1_{\mathbb R_q}(\eta)d\eta
\leq\frac{1}{(2\pi)^n}\int_{\mathbb R_q\cap[-\delta,\delta]}|\det M_{\Phi_{\eta}}|d\eta\\
&\leq\frac{1}{(2\pi)^n}\int_{\mathbb R_{x,q}\cap[-\delta,\delta]}|\det(\mathcal R^L_p+2s\mathcal L_p)|ds.
\end{split}
\end{equation}
Thus we get the conclusion of Theorem \ref{eee}.
\end{proof}

\section{Strong Morse inequalities on CR manifolds with $S^1$ action}

In this section, we will establish the strong Morse inequalities on compact CR manifolds with $S^1$ action. Following the same argument as in Proposition 3.8, Proposition 3.9 in \cite{HM12} and by some minor change we have
\begin{proposition}\label{aaaa}
There exists $u\in\Omega^{0,q}(H_n)$ such that
\begin{equation}
\Box^{(q)}_{b, H_n}u=0,\|u\|_{\Phi_0}=1,
\end{equation}
\begin{equation}
|u(0,0)|^2=(2\pi)^{-n}\int_{\mathbb R_q\cap[-\delta,\delta]}|\det M_{\Phi_\eta}|d\eta,
\end{equation}
and
\begin{equation}
\hat u(z,\eta)\equiv0~\text{when}~|\eta|>\delta.
\end{equation}
\end{proposition}
\begin{proof}
Since some notations have been changed from Proposition 3.8 and 3.9 in \cite{HM12}, we will outline the proof here for the convenient of readers. For any $\eta\in\mathbb R$, we can find a unitary matrix $(a_{ij}(\eta))_{1\leq i,j\leq n-1}$ such that $z_i(\eta)=\sum_{i,j=1}^{n-1}a_{ij}(\eta)z_j$ and  $\Phi_\eta(z)=\sum_{j=1}^{n-1}v_j(\eta)|z_j(\eta)|^2,$ where $v_j(\eta), j=1, \cdots, n-1$ are the eigenvalues of $M_{\Phi_{\eta}}.$ If $\eta\in\mathbb R_q$, we assume $v_1(\eta)<0, \cdots, v_q(\eta)<0, v_{q+1}(\eta)>0, \cdots, v_{n-1}(\eta)>0$. Put
\begin{equation}
\alpha(z,\eta)=C_0|\det M_{\Phi_\eta}|1_{\mathbb R_q\cap[-\delta,\delta]}(\eta)\exp{\left(\sum_{j=1}^qv_j(\eta)|z_j(\eta)|^2\right)}
d\overline {z_1(\eta)}\wedge\cdots d\overline{z_q(\eta)},
\end{equation}
where $C_0=(2\pi)^{1-\frac{n}{2}}\left(\int_{\mathbb R_q\cap[-\delta,\delta]}|\det M_{\Phi_\eta}|d\eta\right)^{-\frac{1}{2}}.$ Then $\Box^{(q)}_{\Phi_{\eta}}\alpha(z,\eta)=0.$ Moreover, we have
\begin{equation}
\int_{\mathbb C^{n-1}}|\alpha(z,\eta)|^2e^{-\Phi_\eta(z)}dv(z)=2\pi\left(\int_{\mathbb R_q\cap[-\delta,\delta]}|\det M_{\Phi_\eta}|d\eta\right)^{-1}|\det M_{\Phi_\eta}|\cdot 1_{\mathbb R_q\cap[-\delta,\delta]}(\eta).
\end{equation}
Set
\begin{equation}
u(z,\theta)=\frac{1}{2\pi}\int_{\mathbb R}\exp{\Big(i\theta\eta+\eta\sum_{j=1}^{n-1}\lambda_j|z_j|^2\Big)}\alpha(z,\eta)d\eta\in\Omega^{0,q}(H_n).
\end{equation}
Using Lemma 3.2 in \cite{HM12}, we can check that $u(z,\theta)$ satifies the properties in Proposition \ref{aaaa}.
\end{proof}
We will use the same notation as in
section \ref{kk}. Fix $x\in X_{{\rm reg\,}}$, choose canonical local coordinates $(z,\theta)$ near $x$ such that $x\leftrightarrow0$, $D=\{(z,\theta)\in\mathbb C^{n-1}\times\mathbb R:|z|<1, |\theta|<\pi\}$. It should be noticed that since $x\in X_{{\rm reg\,}}$, $\theta$ can be defined on $\abs{\theta}<\pi$. Choose two cut-off
functions $\chi\in C_0^\infty(\mathbb C^{n-1}), \tau\in C_0^\infty(\mathbb R)$ in such that $\chi(z)\equiv1$ when $\abs{z}\leq\frac{1}{2}$, $\chi(z)=0$ when $\abs{z}>1$ and $\tau(\theta)\equiv1, |\theta|\leq \frac12; \tau(\theta)\equiv0, |\theta|>1.$ Let $u$ be given as in Proposition \ref{aaaa}. Put
$\chi_{k}(z)=\chi(\frac{z}{\log k}), \tau_k(\theta)=\tau(\frac{\theta}{\log k})$. Put $u_k=\chi_k(\sqrt k z)\tau_k( k\theta)\sum\nolimits^\prime_{|J|=q}u_J(\sqrt kz, k\theta)e^J(z)\in\Omega^{0,q}(X).$ Then ${\rm supp\,}u_k\Subset D_{\frac{\log k}{\sqrt k}}$. Write $\alpha_k=k^{\frac{n}{2}}u_k(z,\theta)\otimes s^k\in\Omega^{0,q}(X, L^k)$. Then
\begin{equation}
\begin{split}
\|\alpha_k\|^2_{h^{L^k}}&=k^n\int_{X}e^{-k\Phi(z)}|u_k(z,\theta)|^2dv_X=k^n\int_Xe^{-k\Phi(z)}
\chi_k^2(\sqrt kz)\tau^2_k(k\theta)|u(\sqrt kz, k\theta)|^2dv_X\\
&=k^n\int_{|z|\leq\frac{\log k}{\sqrt k}, |\theta|\leq\frac{\log k}{k}}e^{-k\Phi(z)}\chi_k^2(\sqrt kz)\tau_k^2(k\theta)|u(\sqrt kz, k\theta)|^2m(z)dv(z)d\theta\\
&=\int_{D_{\log k}}e^{-k\Phi(\frac{z}{\sqrt k})}\chi_k^2(z)\tau_k^2(\theta)|u(z,\theta)|^2m(\frac{z}{\sqrt k})dv(z)d\theta,
\end{split}
\end{equation}
where $m(z)dv(z)d\theta=dv_X$ on $D$.
Then
\begin{equation}
\lim_{k\rightarrow\infty}\|\alpha_k\|^2_{h^{L^k}}=\int_{H_n}e^{-\Phi_0(z)}
|u(z,\theta)|^2dv(z)d\theta=1.
\end{equation}
Second,
\begin{equation}\label{iiii}
k^{-n}|\alpha_k(0,0)|^2_{h^{L^k}}=|u(0,0)|^2=(2\pi)^{-n}\int_{\mathbb R_{x, q}\cap[-\delta,\delta]}|\det(\mathcal R^L_x+2s\mathcal L_x)|ds.
\end{equation}
Third, from $\Box^{(q)}_{b, H_n}u=0$, it is easy to see that there exists a sequence $\mu_k>0$, independent of $p$ and tending to zero such that
\begin{equation}
\left(\frac{1}{k}\Box^{(q)}_{b, k}\alpha_k\big|\alpha_k\right)_{h^{L^k}}\leq \mu_k.
\end{equation}
Moreover, for every $j\in\mathbb N$,
\begin{equation}\label{e-gue150306f}
\mbox{$\left((\frac{1}{k}\Box^{(q)}_{b, k})^j\alpha_k\big|\alpha_k\right)_{h^{L^k}}\To0$ as $k\To\infty$}.
\end{equation}

\begin{theorem}\label{nn}
Set $\beta_k=F_{\delta,k}\alpha_k:=\sum\limits_{|m|\leq k\delta}\alpha_{k,m}$, $T\alpha_{k,m}=(im)\alpha_{k,m}$. Then we will have
\begin{equation}
\begin{split}
&(1)\beta_k\in\Omega^{0,q}_{\leq k\delta}(X, L^k), \lim_{k\rightarrow\infty}\|\beta_k\|^2_{h^{L^k}}=1,\\
&(2)\lim_{k\rightarrow\infty}k^{-n}|\beta_k(x)|^2_{h^{L^k}}=(2\pi)^{-n}\int_{\mathbb R_{x, q}\cap[-\delta, \delta]}|\det(\mathcal R^L_x+2s\mathcal L_x)|ds,\\
&(3)\left(\frac{1}{k}\Box^{(q)}_{b, k}\beta_k|\beta_k\right)_{h^{L^k}}\leq\mu_k,\\
&(4)\,\mbox{$\left((\frac{1}{k}\Box^{(q)}_{b, k})^j\beta_k\big|\beta_k\right)_{h^{L^k}}\To0$ as $k\To\infty$, $\forall j\in\mathbb N$}.
\end{split}
\end{equation}
\end{theorem}
We postpone the proof of Theorem~\ref{nn} until the end of this section.
\begin{proposition}\label{yy}
Let $v_k>0$ be any sequence with $\lim_{k\rightarrow\infty}\frac{\mu_k}{v_k}=0$. Then
\begin{equation}
\liminf_{k\rightarrow\infty} k^{-n}\Pi^{q}_{\leq k\delta, \leq kv_k}(x)\geq (2\pi)^{-n}\int_{\mathbb R_{x,q}\cap[-\delta,\delta]}|\det(\mathcal R^L_x+2s\mathcal L_x)|ds,\ \ \forall x\in X_{{\rm reg\,}}.
\end{equation}
\end{proposition}

\begin{proof}
We will follow the argument of proposition 5.1 in \cite{HM12} to prove this proposition. Let $\mathcal H^{q}_{b, \leq k\delta, >kv_k}(X, L^k)$ denote the space spanned by the eigenforms of $\Box^{(q)}_{b, k}$ restricting to $\Omega^{0, q}_{\leq k\delta}(X, L^k)$ whose eigenvalues  are $>kv_k.$ Fix $x\in X_{{\rm reg\,}}$ and let $\beta_k$ be defined as in Theorem \ref{nn}. $\beta_k=\beta_{k}^1+\beta_k^2$, where
$\beta_k^1\in\mathcal H^q_{\leq k\delta, \leq kv_k}(X, L^k)$, $\beta_k^2\in\overline{\mathcal H^q_{\leq k\delta, >kv_k}(X, L^k)}$. Here the closure of $\mathcal H^q_{\leq k\delta, >kv_k}(X, L^k)$ is under the $Q_b$-norm defined in Proposition 5.1 \cite{HM12}. Then
\begin{equation}\label{e-gueg}
\|\beta_k^2\|^2_{h^{L^k}}=(\beta_k^2|\beta_k^2)_{h^{L^k}}\leq \frac{1}{kv_k}(\Box^{(q)}_{b,k}\beta_k^2|\beta_k^2)_{h^{L^k}}\leq \frac{1}{kv_k}(\Box^{(q)}_{b,k}\beta_k|\beta_k)_{h^{L^k}}\leq\frac{\mu_k}{v_k}\rightarrow0.
\end{equation}
Since $\lim\limits_{k\rightarrow\infty}\|\beta_k\|_{h^{L^k}}=1$, we get $\lim\limits_{k\rightarrow\infty}\|\beta_k^1\|_{h^{L^k}}=1$. Now, we claim that
\begin{equation}
\lim_{k\rightarrow\infty}k^{-n}|\beta_k^2(x)|^2_{h^{L^k}}=0.
\end{equation}
On $D$ with canonical local coordinates $(z, \theta)$, $x\leftrightarrow 0$ and $\Phi(0)=0$, we write $\beta_k^2=k^{\frac{n}{2}}\alpha_k^2\otimes s^k$, $\alpha_k^2\in\Omega^{0,q}(D)$. Then
\begin{equation}\label{ww}
\lim_{k\rightarrow\infty}k^{-n}|\beta_k^2(0)|^2_{h^{L^k}}=\lim_{k\rightarrow\infty}|\alpha_k^2(0)|^2.
\end{equation}
By Proposition \ref{t}, we have
\begin{equation}\label{vv}
|F_k^\ast\alpha_k^2(0)|^2\leq C_{n,r}\left(\|F_k^\ast\alpha_k^2\|^2_{kF_k^\ast\Phi, D_{2r}}+\|\Box^{(q)}_{b,(k)}F_k^\ast\alpha_k^2\|^2_{kF_k^\ast\Phi, n, D_{2r}}+\left\|\left(\frac{\partial}{\partial\theta}\right)^{n}F_k^\ast\alpha_k^2\right\|_{kF_k^\ast\Phi, D_{2r}}\right).
\end{equation}
From the proof of Theorem \ref{ss}, we see that
\begin{equation}\label{tt}
\left\|\left(\frac{\partial}{\partial\theta}\right)^{n}F_k^\ast\alpha_k^2\right\|_{kF_k^\ast\Phi, D_{2r}}\leq
C\|\beta_k^2\|^2_{h^{L^k}}
\end{equation}
where $C>0$ is a constant which does not depend on $k$. Moreover, from (4) in Theorem~\ref{nn},
\begin{equation}\label{uu}
\|\Box^{(q)}_{b,k}F_k^\ast\alpha_k^2\|^2_{kF_k^\ast, n,  D_{2r}}\leq C_1\sum_{m=1}^{n+1}\left\|\left(\frac{1}{k}\Box^{(q)}_{b,k}\right)^m\beta_k^2\right\|^2_{h^{L^k}}\leq C_1\sum_{m=1}^{n+1}\left\|\left(\frac{1}{k}\Box^{(q)}_{b,k}\right)^m\beta_k\right\|^2_{h^{L^k}}\rightarrow0,
\end{equation}
where $C_1>0$ is a constant independent of $k$.
Combining \eqref{e-gueg}, (\ref{ww}), (\ref{vv}), (\ref{tt}) and (\ref{uu}), we get that
\begin{equation}
\lim_{k\rightarrow\infty}k^{-n}|\beta_k^2(0)|^2_{h^{L^k}}=\lim_{k\rightarrow\infty}|\alpha_k^2(0)|^2=0.
\end{equation}
Then
\begin{equation}
\lim_{k\rightarrow\infty} k^{-n}|\beta_{k}^1(0)|^2_{h^{L^k}}=\lim_{k\rightarrow\infty}k^{-n}|\beta_k(0)|^2_{h^{L^k}}=(2\pi)^{-n}
\int_{\mathbb R_{x,q}\cap[-\delta,\delta]}|\det(\mathcal R^L_x+2s\mathcal L_x)|ds.
\end{equation}
Now,
\begin{equation}
k^{-n}\Pi^{q}_{\leq k\delta, \leq kv_k}(x)\geq k^{-n}\frac{|\beta_k^1(x)|^2_{h^{L^k}}}{\|\beta_k^1\|^2_{h^{L^k}}}\rightarrow (2\pi)^{-n}\int_{\mathbb R_{x,q}\cap[-\delta,\delta]}|\det(\mathcal R^L_x+2s\mathcal L_x)|ds.
\end{equation}
The Proposition follows.
\end{proof}

From a simple modification of the proofs of Theorem \ref{ss} and Theorem~\ref{m}, we get the following proposition

\begin{proposition}\label{zz}
Let $v_k>0$ be any sequence with $v_k\rightarrow0$, as $k\rightarrow\infty$. Then there is a constant $C^\prime_0>0$ independent of $k$ such that $k^{-n}\Pi^q_{\leq k\delta,\leq kv_k}(x)\leq C^\prime_0$, $\forall x\in X$, and
\begin{equation}
\limsup_{k\rightarrow\infty} k^{-n}\Pi^{q}_{\leq k\delta, \leq kv_k}(x)\leq (2\pi)^{-n}\int_{\mathbb R_{x,q}\cap[-\delta,\delta]}|\det(\mathcal R^L_x+2s\mathcal L_x)|ds.
\end{equation}
\end{proposition}
Combining Proposition \ref{yy} and Proposition \ref{zz}, we get the conclusion of Theorem \ref{o}.

\begin{proof}[Proof of Theorem~\ref{nn}]
We will use the same notations as in Theorem~\ref{nn}. Note that
$\|\beta_k\|^2_{h^{L^k}}=\sum\limits_{|m|\leq k\delta}\|\alpha_{k,m}\|^2_{h^{L^k}}\leq\|\alpha_k\|^2_{h^{L^k}}.$ On canonical local coordinates $D=\{(z,\theta): |z_j|<r, |\theta|<\pi, j=1,\cdots, n-1\},$ $\alpha_{k,m}$ can be expressed as following:
\begin{equation}\label{bbbb}
\begin{split}
\alpha_{k,m}(z,\theta)&=\frac{1}{2\pi}\int_{-\pi}^{\pi}\alpha_k(z,t)e^{-imt}dte^{im\theta}
=\frac{1}{2\pi}k^{\frac{n}{2}}\int_{-\pi}^{\pi}u_k(z,t)e^{-imt}dte^{im\theta}\otimes s^k\\
&=\frac{1}{2\pi}k^{\frac{n}{2}}\sideset{}{'}\sum_{|J|=q}\int_{-\pi}^{\pi}\chi_k(\sqrt kz)\tau_k(kt)u_J(\sqrt kz, kt)e^{-imt}dte^{im\theta}e^J(z)\otimes s^k\\
&=\frac{1}{4\pi^2}k^{\frac{n}{2}}\int_{-\pi}^{\pi}\int_{|\eta|\leq\delta}\chi_k(\sqrt kz)\hat u(\sqrt kz,\eta)\tau_k(kt)e^{-i(m-k\eta)t}dtd\eta e^{im\theta}\otimes s^k\\
&=\frac{1}{4\pi^2}k^{\frac{n}{2}}\int_{|\eta|\leq\delta}\chi_k(\sqrt kz)\hat u(\sqrt kz,\eta)\hat\tau\left((m-k\eta)\frac{\log k}{k}\right)\frac{\log k}{k}d\eta e^{im\theta}\otimes s^k.
\end{split}
\end{equation}
Assume that $|m|>k\delta$. Then $|m-k\eta|\neq0,$ for all $|\eta|\leq \delta.$ There exists a constant $C>0$ independent of $k$ such that
\begin{equation}\label{cccc}
\left|\hat\tau\left((m-k\eta)\frac{\log k}{k}\right)\right|\leq \frac{C}{|m-k\eta|^3}\frac{k^3}{(\log k)^3}.
\end{equation}
By H\"older inequality
\begin{equation}\label{e-gue150306abI}
\int_{|\eta|\leq\delta}|\hat u(\sqrt kz,\eta)|\frac{1}{|m-k\eta|^3}d\eta\leq\left(\int_{|\eta|\leq\delta}|\hat u(\sqrt kz,\eta)|^2d\eta\right)^{\frac12}\cdot \left(\int_{|\eta|\leq\delta}\frac{1}{(m-k\eta)^6}d\eta\right)^{\frac12}.
\end{equation}
Note that
\begin{equation}\label{e-gue150306abm}
\left(\int_{|\eta|\leq\delta}\frac{1}{(m-k\eta)^6}d\eta\right)^{\frac12}=\frac{1}{\sqrt{5k}}
\left(\left[\frac{1}{(m-k\delta)^5}-\frac{1}{(m+k\delta)^5}\right]\right)^{\frac12}.
\end{equation}
From \eqref{bbbb}, \eqref{cccc}, \eqref{e-gue150306abI} and \eqref{e-gue150306abm}, we have
\begin{equation}
\begin{split}
|\alpha_{k,m}(z,\theta)|_{h^{L^k}}\leq\frac{C}{4\pi^2}k^{\frac{n}{2}}&\chi_k(\sqrt kz)\frac{k^2}{(\log k)^2\sqrt{5k}}e^{-\frac{k\Phi(z)}{2}}\left(\int_{|\eta|\leq \delta}|\hat u(\sqrt kz,\eta)|^2d\eta\right)^{\frac12}\times\\
&\left[\frac{1}{(m-k\delta)^5}-\frac{1}{(m+k\delta)^5}\right]^{\frac12}.
\end{split}
\end{equation}
Let $0<\varepsilon<\frac{1}{3}$ be a small constant. Since
\begin{equation}
\begin{split}
&\sum\limits_{|m|>k\delta+\frac{k}{(\log k)^{1+\varepsilon}}+1}\left[\frac{1}{(m-k\delta)^5}-\frac{1}{(m+k\delta)^5}\right]^{\frac12}\\
&=2\sum\limits_{m>k\delta+\frac{k}{(\log k)^{1+\varepsilon}}+1}\left[\frac{1}{(m-k\delta)^5}-\frac{1}{(m+k\delta)^5}\right]^{\frac12}\\
&\leq 2\sum\limits_{m>k\delta+\frac{k}{(\log k)^{1+\varepsilon}}+1}\left[\frac{1}{(m-k\delta)}\right]^{\frac52}\\
&\leq2\int_{\frac{k}{(\log k)^{1+\varepsilon}}}^\infty\frac{1}{\sigma^{\frac52}}d\sigma=\frac{4}{3}\frac{(\log k)^{\frac32(1+\varepsilon)}}{k^{\frac32}},
\end{split}
\end{equation}
we have
\begin{equation}
\begin{split}
&\Big|\sum\limits_{|m|>k\delta+\frac{k}{(\log k)^{1+\varepsilon}}+1}\alpha_{k,m}(z,\theta)\Big|_{h^{L^k}}\\
&\leq\frac{C}{3\pi^2}k^{\frac{n}{2}}\chi_k(\sqrt kz)\frac{k^2}{(\log k)^2}\frac{1}{\sqrt{5k}}\frac{(\log k)^{\frac32(1+\varepsilon)}}{k^{\frac32}}
e^{-\frac{k\Phi(z)}{2}}\left(\int_{|\eta|\leq\delta}|\hat u(\sqrt kz,\eta)|^2d\eta\right)^{\frac12}\\
&\leq C_1k^{\frac{n}{2}}\chi_k(\sqrt kz)\frac{1}{(\log k)^{\frac{1-3\varepsilon}{2}}}e^{-\frac{k\Phi(z)}{2}}\left(\int_{|\eta|\leq\delta}|\hat u(\sqrt kz,\eta)|^2d\eta\right)^{\frac12},
\end{split}
\end{equation}
where $C_1>0$ is a constant independent of $k$.
Let $\gamma_k(z,\theta)=\sum\limits_{|m|>k\delta+\frac{k}{(\log k)^{1+\varepsilon}}+1}\alpha_{k,m}(z,\theta)$. Then
\begin{equation}\label{eeee}
|\gamma_k(z,\theta)|^2_{h^{L^k}}\leq C_1k^n\chi_k^2(\sqrt kz)\frac{1}{(\log k)^{1-3\varepsilon}}e^{-k\Phi(z)}\int_{|\eta|\leq\delta}|\hat u(\sqrt k z,\eta)|^2d\eta.
\end{equation}
For any $M>0$,
\begin{equation}
\begin{split}
&\int_{|z|\leq\frac{\log k}{k}}\int_{|\theta|\leq\frac{M}{k}}|\gamma_k(z,\theta)|^2_{h^{L^k}}m(z)dv(z)d\theta\\
&\leq C_2\frac{M}{(\log k)^{1-3\varepsilon}}\int_{|z|\leq \log k}\int_{|\eta|\leq\delta}\chi_k^2(z)e^{-k\Phi(\frac{z}{\sqrt k})}|\hat u(z,\eta)|^2m(\frac{z}{\sqrt k})dv(z)d\theta.
\end{split}
\end{equation}
Letting $k\rightarrow\infty$, we have
\begin{equation}\label{e-guebpa}
\lim_{k\rightarrow\infty}\int_{|z|\leq\frac{\log k}{\sqrt k}}\int_{|\theta|\leq \frac{M}{k}}|\gamma_k(z,\theta)|^2_{h^{L^k}}m(z)dv(z)d\theta=0.
\end{equation}
On the other hand,
\begin{equation}
\begin{split}
&|\alpha_{k,m}(z,\theta)|_{h^{L^k}}\leq\frac{1}{2\pi}k^{\frac{n}{2}}\int_{-\pi}^{\pi}\int_{|\eta|\leq\delta}
\chi_k(\sqrt kz)\tau_k(kt)|\hat u(\sqrt kz,\eta)|dtd\eta e^{-\frac{k\Phi(z)}{2}}\\
&\leq\frac{1}{2\pi}k^{\frac{n}{2}}\int_{|\eta|\leq\delta}|\hat u(\sqrt kz,\eta)|d\eta\chi_k(\sqrt kz)\left(\frac{\log k}{k}\right)d\eta e^{-\frac{k\Phi(z)}{2}}\\
&\leq Ck^{\frac{n}{2}}\left(\frac{\log k}{k}\right)\left(\int_{|\eta|\leq\delta}|\hat u(\sqrt kz,\eta)|^2d\eta\right)^{\frac12}d\eta\chi_k(\sqrt kz)e^{\frac{-k\Phi(z)}{2}},
\end{split}
\end{equation}
where $C>0$ is a constant independent of $k$.
Let $\sigma_k=\sum\limits_{k\delta<|m|\leq k\delta+\frac{k}{(\log k)^{1+\varepsilon}}}\alpha_{k,m}(z,\theta)$. Then
\begin{equation}\label{ffff}
\begin{split}
|\sigma_k|_{h^{L^k}}
&\leq C \frac{\log k}{k}\frac{k}{(\log k)^{1+\varepsilon}}k^{\frac{n}{2}}\left(\int_{|\eta|\leq\delta}|\hat u(\sqrt kz,\eta)|^2d\eta\right)^{\frac12}\chi_k(\sqrt kz)e^{\frac{-k\Phi(z)}{2}}\\
&\leq C\frac{1}{(\log k)^\varepsilon}k^{\frac{n}{2}}\left(\int_{|\eta|\leq\delta}|\hat u(\sqrt kz,\eta)|^2d\eta\right)^{\frac12}\chi_k(\sqrt kz)e^{-\frac{k\Phi(z)}{2}}.
\end{split}
\end{equation}
From \eqref{ffff}, we can check that
\begin{equation}\label{e-guebpII}
\lim_{k\rightarrow\infty}\int_{|z|\leq\frac{\log k}{\sqrt k}, |\theta|\leq \frac{M}{k}}|\sigma_k|^2m(z)dv(z)=0,\ \ \forall M>0.
\end{equation}
Write $\alpha_k=\beta_k+\gamma_k+\sigma_k$. Here, $\beta_k=\sum\limits_{|m|\leq k\delta}\alpha_{k,m}$. Then
\begin{equation}\label{e-guebpIII}
\int_{|z|\leq\frac{\log k}{\sqrt k}, |\theta|\leq\frac{M}{k}}|\beta_k|^2_{h^{L^k}}m(z)dv(z)d\theta=\int_{|z|\leq\frac{\log k}{\sqrt k}, |\theta|\leq\frac{M}{k}}|\alpha_k-\gamma_k-\sigma_k|^2_{h^{L^k}}m(z)dv(z)d\theta.
\end{equation}
Since
\begin{equation}\label{e-guebp}
\lim_{k\rightarrow\infty}\int_{|z|\leq\frac{\log k}{\sqrt k}, |\theta|\leq\frac{M}{k}}|\alpha_k|^2_{h^{L^k}}m(z)dv(z)d\theta=\int_{\mathbb C^{n-1}\times\{\theta\in\mathbb R: |\theta|\leq M\}}|u(z,\theta)|^2e^{-\Phi_0(z)}dv(z)d\theta
\end{equation}
and $\|u\|^2=1$, for any $\epsilon>0$, we can choose a constant $M>0$ such that
\begin{equation}\label{e-guebpI}
\int_{\mathbb C^{n-1}\times\{\theta\in\mathbb R: |\theta|\leq M\}}|u(z,\theta)|^2e^{-\Phi_0(z)}dv(z)d\theta\geq1-\epsilon.
\end{equation}
From \eqref{e-guebpa}, \eqref{e-guebpII}, \eqref{e-guebpIII}, \eqref{e-guebp} and \eqref{e-guebpI}, we deduce that
\begin{equation}
\liminf_{k\rightarrow\infty} \|\beta_k\|^2_{h^{L^k}}\geq\lim_{k\rightarrow\infty}\int_{|z|\leq\frac{\log k}{\sqrt k}, |\theta|\leq\frac{M}{k}}|\beta_k|^2_{h^{L^k}}m(z)dv(z)d\theta\geq 1-\epsilon,\ \ \forall\epsilon>0.
\end{equation}
Thus, $$\liminf_{k\rightarrow\infty}\|\beta_k\|^2_{h^{L^k}}\geq 1.$$ On the other hand, $\|\beta_k\|^2_{h^{L^k}}\leq\|\alpha_k\|^2_{h^{L^k}}\leq 1$, then we have
\begin{equation}
\lim_{k\rightarrow\infty}\|\beta_k\|^2_{h^{L^k}}=1.
\end{equation}

Proof of (2) in Theorem \ref{nn}. Recall that $\alpha_k=\beta_k+\gamma_k+\sigma_k$. From (\ref{eeee}), we have
\begin{equation}
k^{-n}|\gamma_k(0,0)|^2_{h^{L^k}}\leq C_1\frac{1}{(\log k)^{1-3\varepsilon}}\int_{|\eta|\leq\delta}|\hat u(0,\eta)|^2d\eta.
\end{equation}
Then
\begin{equation}\label{gggg}
\lim_{k\rightarrow\infty}k^{-n}|\gamma_k(0,0)|^2_{h^{L^k}}=0.
\end{equation}
From (\ref{ffff}), we have
\begin{equation}
k^{-n}|\sigma_k(0,0)|^2_{h^{L^k}}\leq
\int_{|\eta|\leq\delta}|\hat u(0,\eta)|^2d\eta\frac{1}{(\log k)^{2\varepsilon}}.
\end{equation}
Then
\begin{equation}\label{hhhh}
\lim_{k\rightarrow\infty}k^{-n}|\sigma_k(0,0)|^2_{h^{L^k}}=0.
\end{equation}
Combining (\ref{iiii}), (\ref{gggg}) and (\ref{hhhh}) we get the conclusion of the second part of Theorem \ref{nn}.

Proof of (3) in Theorem \ref{nn}. $\beta_k=F_{\delta,k}\alpha_k$. Since $\alpha_k=F_{\delta,k}\alpha_k+(I-F_{\delta, k})\alpha_k.$ Since $\Box^{(q)}_{b, k}F^{(q)}_{\delta, k}=F^{(q)}_{\delta, k}\Box^{(q)}_{b, k}$, then
\begin{equation}
\Box^{(q)}_{b, k}\alpha_k=\Box^{(q)}_{b, k}F_{\delta ,k}\alpha_k+(I-F_{\delta, k})\Box^{(q)}_{b,k}\alpha_k
\end{equation}
and
\begin{equation}
\begin{split}
\left(\frac{1}{k}\Box^{(q)}_{b,k}\beta_k|\beta_k\right)_{h^{L^k}}
&=\left(\frac{1}{k}\Box^{(q)}_{b,k}\alpha_k|F_{\delta,k}\alpha_k\right)_{h^{L^k}}
=\left(F_{\delta, k}\frac{1}{k}\Box^{(q)}_{b,k}\alpha_k|F_{\delta, k}\alpha_k\right)_{h^{L^k}}\\
&\leq\left(\frac{1}{k}\Box_{b,k}^{(q)}\alpha_k|\alpha_k\right)_{h^{L^k}}\leq\mu_k.
\end{split}
\end{equation}
for some $\mu_k$ tending to zero. Similarly, we can repeat the procedure above and get (4) in Theorem~\ref{nn}. Thus, we get the conclusion of Theorem of \ref{nn}.
\end{proof}

\section{Examples}

In this section, some examples are collected. The aim is to illustrate the main results in some
simple situations.

\subsection{CR manifolds in projective spaces}\label{s-cmips}

We consider $\Complex\mathbb P^{N-1}$, $N\geq3$. Let $[z]=[z_1,\ldots,z_N]$ be the homogeneous coordinates of $\Complex\mathbb P^{N-1}$. Put
\[X:=\set{[z_1,\ldots,z_N]\in\Complex\mathbb P^{N-1};\, \lambda_1\abs{z_1}^2+\cdots+\lambda_m\abs{z_m}^2+\cdots+\lambda_N\abs{z_N}^2=0},\]
where $m\in\mathbb N$ and $\lambda_j\in\Real$ $j=1,\ldots,N$.  We assume that $\lambda_1<0,\ldots,\lambda_m<0$, $\lambda_{m+1}>0,\lambda_{m+2}>0,\ldots,\lambda_N>0$. Then $X$ is a compact CR manifold of dimension $2(N-1)-1$ with CR structure $T^{1,0}X:=T^{1,0}\Complex\mathbb P^{N-1}\bigcap\Complex TX$. $X$ admits a $S^1$ action:
\begin{equation}\label{e-gue131218III}
\begin{split}
S^1\times X&\To X,\\
e^{i\theta}\circ[z_1,\ldots,z_m,z_{m+1},\ldots,z_N]&\To[e^{i\theta}z_1,\ldots,e^{i\theta}z_m,z_{m+1},\ldots,z_N],\ \ \theta\in[-\pi,\pi).
\end{split}
\end{equation}
Since $(z_1,\ldots,z_m)\neq0$ on $X$, this $S^1$ action is well-defined. Moreover, it is straightforward to check that this $S^1$ action is CR and transversal. Let $T$ be the global vector field induced by the $S^1$ action.

Let $E\To\Complex\mathbb P^{N-1}$ be the hyperplane line bundle with respect to the Fubini-Study metric. For $j=1,2,\ldots,N$, put $W_j=\set{[z_1,\ldots,z_N]\in\Complex\mathbb P^{N-1};\, z_j\neq0}$. Then, $E$ is trivial on $W_j$, $j=1,\ldots,N$, and we can find local trivializing section $e_j$ of $E$ on $W_j$, $j=1,\ldots,N$, such that for every $j, t=1,\ldots,N$,
\begin{equation}\label{e-gue131218III-I}
e_j(z)=\frac{z_j}{z_t}e_t(z)\ \ \mbox{on $W_j\bigcap W_t$},\ \ z=[z_1,\ldots,z_N]\in W_j\bigcap W_t.
\end{equation}
Consider $L:=E|_X$. Then, $L$ is a CR line bundle over $(X,T^{1,0}X)$. It is easy to see that $X$ can be covered with open sets $U_j:=W_j|_X$, $j=1,2,\ldots,m$, with trivializing sections $s_j:=e_j|_X$, $j=1,2,\ldots,m$, such that the corresponding transition functions are rigid CR functions. Thus, $L$ is a rigid CR line bundle over $(X,T^{1,0}X)$. Let $h^L$ be the Hermitian fiber metric on $L$ given by
\[\abs{s_j(z_1,\ldots,z_N)}^2_{h^L}:=e^{-\log\bigr(\frac{\abs{z_1}^2+\cdots+\abs{z_N}^2}{\abs{z_j}^2}\bigr)},\ \ j=1,\ldots,m.\]
It is not difficult to check that $h^L$ is well-defined and $h^L$ is a rigid positive CR line bundle.  From this and Theorem~\ref{llll}, we conclude that $L$ is a big line bundle over $X$.

\subsection{Compact Heisenberg groups} \label{s-chg}

Let $\lambda_1,\ldots,\lambda_{n-1}$ be given non-zero integers.
Let $\mathscr CH_n=(\Complex^{n-1}\times\Real)/_\sim$\,, where
$(z, t)\sim(\Td z, \Td t)$ if
\[\begin{split}
&\Td z-z=(\alpha_1,\ldots,\alpha_{n-1})\in\sqrt{2\pi}\mathbb Z^{n-1}+i\sqrt{2\pi}\mathbb Z^{n-1},\\
&\Td t-t-i\sum^{n-1}_{j=1}\lambda_j(z_j\ol\alpha_j-\ol z_j\alpha_j)\in 2\pi\mathbb Z.
\end{split}\]
We can check that $\sim$ is an equivalence relation
and $\mathscr CH_n$ is a compact manifold of dimension $2n-1$. The equivalence class of $(z,t)\in\Complex^{n-1}\times\Real$ is denoted by
$[(z, t)]$. For a given point $p=[(z, t)]$, we define
$T^{1,0}_p\mathscr CH_n$ to be the space spanned by
\[
\textstyle
\big\{\frac{\pr}{\pr z_j}+i\lambda_j\ol z_j\frac{\pr}{\pr t},\ \ j=1,\ldots,n-1\big\}.
\]
It is easy to see that the definition above is independent of the choice of a representative $(z,t)$ for $[(z,t)]$.
Moreover, we can check that $T^{1,0}\mathscr CH_n$ is a CR structure. $\mathscr CH_n$ admits the natural $S^1$ action: $e^{i\theta}\circ [z,t]\To [z,t+\theta]$, $0\leq\theta<2\pi$. Let $T$ be the global vector field induced by this $S^1$ action. We can check that this $S^1$ action is CR and transversal and
$T=\frac{\pr}{\pr t}$.  We take a Hermitian metric $\langle\,\cdot\,|\,\cdot\,\rangle$ on the complexified tangent bundle $\Complex T\mathscr CH_n$ such that
\[
\Big\lbrace
\tfrac{\pr}{\pr z_j}+i\lambda_j\ol z_j\tfrac{\pr}{\pr t}\,, \tfrac{\pr}{\pr\ol z_j}-i\lambda_jz_j\tfrac{\pr}{\pr t}\,, -\tfrac{\pr}{\pr t}\,;\, j=1,\ldots,n-1\Big\rbrace
\]
 is an orthonormal basis. The dual basis of the complexified cotangent bundle is
\[
\Big\lbrace
dz_j\,,\, d\ol z_j\,,\, \omega_0:=-dt+\textstyle\sum^{n-1}_{j=1}(i\lambda_j\ol z_jdz_j-i\lambda_jz_jd\ol z_j); j=1,\ldots,n-1
\Big\rbrace\,.
\]
The Levi form $\mathcal{L}_p$ of $\mathscr CH_n$ at $p\in\mathscr CH_n$ is given by
$\mathcal{L}_p=\sum^{n-1}_{j=1}\lambda_jdz_j\wedge d\ol z_j$.

Now, we construct a rigid CR line bundle $L$ over $\mathscr CH_n$. Let $L=(\Complex^{n-1}\times\Real\times\Complex)/_\equiv$ where $(z,\theta,\eta)\equiv(\Td z, \Td\theta, \Td\eta)$ if
\[\begin{split}
&(z,\theta)\sim(\Td z, \Td\theta),\\
&\Td\eta=\eta\exp(\sum^{n-1}_{j,t=1}\mu_{j,t}(z_j\ol\alpha_t+\frac{1}{2}\alpha_j\ol\alpha_t)),
\end{split}\]
where $\alpha=(\alpha_1,\ldots,\alpha_{n-1})=\Td z-z$, $\mu_{j,t}=\mu_{t,j}$, $j, t=1,\ldots,n-1$, are given integers. We can check that $\equiv$ is an equivalence relation and
$L$ is a rigid CR line bundle over $\mathscr CH_n$. For $(z, \theta, \eta)\in\Complex^{n-1}\times\Real\times\Complex$, we denote
$[(z, \theta, \eta)]$ its equivalence class.
It is straightforward to see that the pointwise norm
\[
\big\lvert[(z, \theta, \eta)]\big\rvert^2_{h^L}:=\abs{\eta}^2\exp\big(-\textstyle\sum^{n-1}_{j,t=1}\mu_{j,t}z_j\ol z_t\big)
\]
is well defined. In local coordinates $(z, \theta, \eta)$, the weight function of this metric is
\[\phi=\sum^{n-1}_{j,t=1}\mu_{j,t}z_j\ol z_t.\]
Thus, $L$ is a rigid CR line bundle over $\mathscr CH_n$ with rigid Hermitian metric $h^L$.
Note that
\[
\textstyle\ddbar_b=\sum^{n-1}_{j=1}d\ol z_j\wedge(\frac{\pr}{\pr\ol z_j}-i\lambda_jz_j\frac{\pr}{\pr\theta})\,,\quad
\pr_b=\sum^{n-1}_{j=1}dz_j\wedge(\frac{\pr}{\pr z_j}+i\lambda_j\ol z_j\frac{\pr}{\pr\theta}).
\]
Thus
$d(\ddbar_b\phi-\pr_b\phi)=2\sum^{n-1}_{j,t=1}\mu_{j,t}dz_j\wedge d\ol z_t$ and for any $p\in\mathscr CH_n$,
\[\mathcal R^L_p=\sum^{n-1}_{j,t=1}\mu_{j,t}dz_j\wedge d\ol z_t.\]
 From this and Theorem~\ref{llll}, we conclude that

\begin{theorem}\label{t-gue131221}
If $\left(\mu_{j,t}\right)^{n-1}_{j,t=1}$ is positive definite, then $L$ is a big line bundle on $\mathscr CH_n$.
\end{theorem}

\subsection{Holomorphic line bundles over a complex torus}\label{s-hlbo}

Let
\[T_n:=\Complex^n/(\sqrt{2\pi}\mathbb Z^n+i\sqrt{2\pi}\mathbb Z^n)\]
be the flat torus. Let $\lambda=\left(\lambda_{j,t}\right)^{n}_{j,t=1}$, where $\lambda_{j,t}=\lambda_{t,j}$,
$j, t=1,\ldots,n$, are given integers. Let $L_\lambda$ be the holomorphic
line bundle over $T_n$
with curvature the $(1,1)$-form
$\Theta_\lambda=\sum^n_{j,t=1}\lambda_{j,t}dz_j\wedge d\ol z_t$.
More precisely, $L_\lambda:=(\Complex^n\times\Complex)/_\sim$\,, where
$(z, \theta)\sim(\Td z, \Td\theta)$ if
\[
\Td z-z=(\alpha_1,\ldots,\alpha_n)\in \sqrt{2\pi}\mathbb Z^n+i\sqrt{2\pi}\mathbb Z^n\,,\quad
\Td\theta=\textstyle\exp\big(\sum^n_{j,t=1}\lambda_{j,t}(z_j\ol\alpha_t+\tfrac{1}{2}\alpha_j\ol\alpha_t\,)\big)\theta\,.
\]
We can check that $\sim$ is an equivalence relation and $L_\lambda$ is a holomorphic line bundle over $T_n$.
For $[(z, \theta)]\in L_\lambda$
we define the Hermitian metric by
\[
\big\vert[(z, \theta)]\big\vert^2:=\abs{\theta}^2\textstyle\exp(-\sum^n_{j,t=1}\lambda_{j,t}z_j\ol z_t)
\]
and it is easy to see that this definition is independent of the choice of a representative $(z, \theta)$ of $[(z, \theta)]$. We denote by $\phi_\lambda(z)$ the weight of this Hermitian fiber metric. Note that $\pr\ddbar\phi_\lambda=\Theta_\lambda$.

Let $L^*_\lambda$ be the
dual bundle of $L_\lambda$ and let $\norm{\,\cdot\,}_{L^*_\lambda}$ be the norm of $L^*_\lambda$ induced by the Hermitian fiber metric on $L_\lambda$. Consider the compact CR manifold of dimension $2n+1$: $X=\{v\in L^*_\lambda:\, \norm{v}_{L^*_\lambda}=1\}$; this is the boundary of the Grauert tube associated to $L^*_\lambda$. The manifold $X$ is equipped with a natural $S^1$-action.
Locally $X$ can be represented in local holomorphic coordinates $(z,\eta)$, where $\eta$ is the fiber coordinate, as the set of all $(z,\eta)$ such that $\abs{\eta}^2e^{\phi_\lambda(z)}=1$.
The $S^1$-action on $X$ is given by $e^{i\theta}\circ (z,\eta)=(z,e^{i\theta}\eta)$, $e^{i\theta}\in S^1$, $(z,\eta)\in X$. Let $T$ be the global vector field on $X$ induced by this $S^1$ action. We can check that this $S^1$ action is CR and transversal.

Let $\pi:L^*_\lambda\To T_n$
be the natural projection from $L^*_\lambda$ onto $T_n$. Let $\mu=\left(\mu_{j,t}\right)^{n}_{j,t=1}$, where $\mu_{j,t}=\mu_{t,j}$, $j, t=1,\ldots,n$, are given integers. Let $L_\mu$ be another holomorphic
line bundle over $T_n$ determined by the constant curvature form
$\Theta_\mu=\sum^n_{j,t=1}\mu_{j,t}dz_j\wedge d\ol z_t$ as above.
The pullback line bundle $\pi^*L_\mu$ is a holomorphic line bundle over $L^*_\lambda$. If we restrict $\pi^*L_\mu$ on $X$, then we can check that $\pi^*L_\mu$ is a rigid CR line bundle over $X$.

The Hermitian fiber metric on $L_\mu$ induced by $\phi_\mu$ induces a Hermitian fiber metric on $\pi^*L_\mu$ that we shall denote by $h^{\pi^*L_\mu}$. We let $\psi$ to denote the weight of $h^{\pi^*L_\mu}$.
The part of $X$ that lies over a fundamental domain of $T_n$ can be represented in local holomorphic coordinates
$(z, \xi)$, where $\xi$ is the fiber coordinate, as the set of all $(z, \xi)$ such that
$r(z, \xi):=\abs{\xi}^2\exp(\sum^n_{j,t=1}\lambda_{j,t}z_j\ol z_t)-1=0$
and the weight $\psi$ may be written as $\psi(z, \xi)=\sum^n_{j,t=1}\mu_{j,t}z_j\ol z_t$. For convenient we denote $\pi^\ast L_{\mu}$ by $L$.
From this we see that $L$ is a $T-$rigid CR line bundle over $X$ with rigid Hermitian fiber metric $h^{L}$. It is straightforward to check that for any $p\in X$, we have
$\mathcal R^{L}_p=\frac12 d(\ddbar_b\psi-\pr_b\psi)(p)|_{T^{1, 0}X}=\sum^n_{j,t=1}\mu_{j,t}dz_j\wedge d\ol z_t$.
Thus, if $\left(\mu_{j,t}\right)^{n-1}_{j,t=1}$ is positive definite, then $L$ is a rigid positive CR line bundle. From this and Theorem~\ref{llll}, we conclude that

\begin{theorem}\label{t-gue131221I}
If $\left(\mu_{j,t}\right)^{n-1}_{j,t=1}$ is positive definite, then $L$ is a big line bundle over $X$.
\end{theorem}

\begin{center}
{\bf Acknowledgement}
\end{center}

The authors would like to thank the referee for many detailed remarks that have helped to improve the presentation.

\bibliographystyle{amsalpha}

\end{document}